\documentclass[aip,jcp,groupedaddress,preprint]{revtex4-1}%
\usepackage{amssymb,amsmath,amsthm,amstext,amscd,bm}
\usepackage[mathcal]{euscript}
\usepackage{graphicx}%
\usepackage{mathtools}%
    \mathtoolsset{showonlyrefs}%
\usepackage[normalem]{ulem}
\usepackage{graphicx}%
\usepackage{subfig}%
\usepackage{color}%

\newcommand{\shrinkfac}{0.6} 
\newcommand{\shrinkfactwo}{0.3} 

\DeclareMathOperator{\pexp}{\mathbf{E}}
\DeclareMathOperator{\pvar}{Var}
\DeclareMathOperator{\pcov}{Cov}
\newcommand{\dd}{\mathrm{d}}
\DeclareMathOperator{\trace}{tr}
\DeclareMathOperator{\bias}{Bias}
\DeclareMathOperator{\mse}{MSE}
\DeclareMathOperator{\VACF}{VACF}
\DeclareMathOperator{\PACF}{PACF}
\DeclareMathOperator{\MSD}{MSD}
\DeclareMathOperator{\AR}{AR}

\theoremstyle{plain}%
\newtheorem{opprob}{Optimization Problem} 

\begin{document}

\title{Uncertainty quantification for generalized Langevin dynamics}

\author{Eric J.~\surname{Hall}} \email{hall@math.umass.edu}
\author{Markos A.~\surname{Katsoulakis}} \email{markos@math.umass.edu}
\author{Luc \surname{Rey-Bellet}} \email{luc@math.umass.edu}

\affiliation{Department of Mathematics and Statistics, University of
  Massachusetts Amherst\\ Amherst, Massachusetts 01003, USA}

\begin{abstract} We present efficient finite difference estimators for
  goal-oriented sensitivity indices with applications to the
  generalized Langevin equation (GLE). In particular, we apply these
  estimators to analyze an extended variable formulation of the GLE
  where other well known sensitivity analysis techniques such as the
  likelihood ratio method are not applicable to key parameters of
  interest. These easily implemented estimators are formed by coupling
  the nominal and perturbed dynamics appearing in the finite
  difference through a common driving noise, or common random
  path. After developing a general framework for variance reduction
  via coupling, we demonstrate the optimality of the common random
  path coupling in the sense that it produces a minimal variance
  surrogate for the difference estimator relative to sampling dynamics
  driven by independent paths. In order to build intuition for the
  common random path coupling, we evaluate the efficiency of the
  proposed estimators for a comprehensive set of examples of interest
  in particle dynamics. These reduced variance difference estimators
  are also a useful tool for performing global sensitivity analysis
  and for investigating non-local perturbations of parameters, such as
  increasing the number of Prony modes active in an extended variable
  GLE.
\end{abstract}



\maketitle%

\section{Introduction}%
\label{sec:introduction}%

Sensitivity analysis (SA), understanding how changes in input
parameters affect the output of a system, is a key component of
uncertainty quantification (UQ), optimal experimental design, and
analysis of model robustness, identifiability, and reliability.
\cite{SaltelliEtAl:2000sa,SaltelliEtAl:2000sb} 
The local sensitivity of a system can be analyzed by computing
sensitivity indices that are formed by taking partial derivatives with
respect to each of the
input parameters. These indices quantify which parameter directions
are most sensitive to perturbations.

The present articles concerns SA techniques for the GLE and other
models of interest in particle dynamics. The Langevin equation (LE)
models particle diffusion in the presence of a heat bath where the
particle-bath interactions are reduced to an instantaneous drag force
and a delta-correlated random force.\cite{CoffeyKalmykov:2012le} This
approximation dramatically reduces the computational cost compared to
explicitly resolving the particle-bath interactions. However, there
are a number of compelling applications where the Langevin assumptions
fail to produce a reliable model, such as anomalous diffusion. The
GLE, a more reliable model of anomalous diffusion, incorporates
``memory'' into the drag force through the inclusion of a kernel
depending on the history of the velocity. In many instances, this
non-Markovian system can be mapped onto a Markovian system with
additional degrees of freedom under physically reasonable assumptions,
such as when the memory kernel can be approximated by a positive Prony
series.\cite{BaczewskiBond:2013ni} The resulting extended variable
formulation contains many parameters that must be tuned and is
therefore an ideal candidate for SA and UQ. However, well known SA
techniques such as likelihood ratio and pathwise methods are not
applicable to analyze the sensitivity of key parameters of interest in
the extended variable formulation. In contrast, Monte Carlo finite
difference estimators of sensitivity indices are applicable to all
parameters of interest in the extended variable GLE, but introduce a
bias error and typically have a large variance making them
computationally expensive.

We give efficient Monte Carlo finite difference estimators via a
coupling method for approximating goal-oriented sensitivity indices
for a large class of SDEs. In particular, we apply these estimators to
an extended variable formulation of the GLE where the memory kernel
can be approximated by a positive Prony series, a choice motivated by
applications in anomalous diffusion in biological
fluids.\cite{MasonWeitz:1995ve, Goychuk:2012ve} In the context of this
application area, we mention that other authors have given a Bayesian
methodology for comparing different models of anomalous diffusion that
favors the GLE.\cite{LysyEtAl:2016cm} In addition to biological
fluids, other recent interesting applications of the GLE include:
modeling nanoscale materials and solids; \cite{Kantorovich:2008oi,
  KantorovichRompotis:2008ii, StellaEtAl:2014gl, NessEtAl:2015ap,
  NessEtAl:2016ne} thermostats for sampling classical and path
integral molecular dynamics;\cite{CeriottiEtAl:2011ac,
  MorroneEtAl:2011mt, CeriottiEtAl:2010cn, CeriottiEtAl:2010et,
  CeriottiEtAl:2009lc, CeriottiEtAl:2009qe} and, more generally,
reduced order modeling.\cite{ChorinStinis:2007me, HijonEtAl:2010mz,
  DarveEtAl:2009co, LiEtAl:2015mz} This list of applications is far
from exhaustive but nevertheless provides strong incentive for
investigating SA and UQ techniques for the GLE and its extended
variable formulations.

To provide further orientation consider, for simplicity, the
sensitivity of the stochastic dynamics $X_t(\theta)$, depending on an
input parameter $\theta$,
\begin{equation*}
  \mathcal{S}(t,\theta;f) = \partial_\theta \pexp[f(X_t(\theta))],
\end{equation*}
for a given observable $f$ where $\partial_\theta$ is the derivative
with respect to $\theta$. In general, the finite difference approach
is to approximate the derivative above by a finite difference quotient
and then obtain the required moments by Monte Carlo. For example, a
forward difference with bias parameter $\varepsilon$ yields the
sensitivity estimator,
\begin{equation}
  \label{eq:1}
  \mathcal{S}_\varepsilon(t,\theta;f) 
  = \left(\pexp[f(X_t(\theta + \varepsilon))] 
    - \pexp[f(X_t(\theta))]\right) / \varepsilon,
\end{equation}
where $\mathcal{S}_\varepsilon \approx \mathcal{S}$ for $\varepsilon$
small, and then the estimator is computed by approximating the
expectations with sample averages. Similar expressions can be given
for central differences and more general finite difference
stencils. While this approach requires little analysis of the
underlying model and is easily implemented, the introduction of the
bias, and in particular its effect on the variance of
$\mathcal{S}_\varepsilon$, a key quantity in evaluating the efficiency
of the method, is often cited as a reason for pursuing alternative
methods.\cite{ChenGlasserman:2007ma} However, as we shall show, the
variance of $\mathcal{S}_\varepsilon$ can be reduced by choosing the
right sampling strategy for the observable of the nominal and
perturbed dynamics, respectively, $f(X_t(\theta))$ and
$f(X_t(\theta+\varepsilon))$ in the expression above. For a
comprehensive set of examples of interest in particle dynamics, we
demonstrate that coupling the nominal and perturbed dynamics through a
common driving noise, that is, a \emph{common random path coupling},
reduces the variance of the finite difference estimator, often
substantially. In particular, for the extended variable GLE with a
\emph{convex} potential the reduction due to the common random path
coupling is on the order of the bias squared---mitigating the effect
of the systematic error. The common random path coupling also leads to
reduced variance estimators for problems with \emph{nonconvex}
potentials, although the reduction is not expected to be on the order
of the bias squared (cf. Figures~\ref{fig:sa_c0_varofvacf} and
\ref{fig:gle_nonconvex}). This is a topic that deserves further
rigorous analysis that will be the subject of future work.

Other well known SA techniques for continuous time stochastic
dynamics, including pathwise methods,
\cite{PantazisKatsoulakis:2013re} likelihood ratio methods,
\cite{Glynn:1990lr,ArampatzisEtAl:2015lr} and Malliavin
methods,\cite{FournieEtAl:1999oi,FournieEtAl:2001ii} produce
\emph{unbiased} estimators of sensitivity indices by representing
$\mathcal{S}$ as the expectation of $f$ (or its derivative) under a
change of measure. A very good account of the interrelations among
them has been given by other authors.\cite{ChenGlasserman:2007ma}
However, each of these methods is not suited to our application of
interest, the GLE, for reasons that we detail below.

The pathwise and likelihood ratio methods are not applicable to key
parameters of interest, those appearing in the drift and diffusion
terms, in the extended variable formulation of the GLE.  The pathwise
method views the dynamics at each fixed time as a density and takes a
derivative of this density yielding the estimator,
$\mathcal{S}_{P}(t,\theta;f) = \partial_\theta \pexp[f(X_t(\theta))] =
\pexp \left[ f^\prime(X_t(\theta)) \partial_\theta
  X_t(\theta)\right],$ requiring equality to hold when the order of
differentiation and expectation are interchanged and a smooth
observable $f$. In its most general form, if an expression
$\pexp[f(X_t(\theta))] = \int f(x_t) g(\theta, x_t) \dd x_t$ exists,
then the likelihood ratio estimator,
\begin{align*}
  \mathcal{S}_{LR}(t,\theta;f) 
  &= \partial_\theta \pexp[f(X_t(\theta))]\\ 
  &= \int f(x_t) \left[ \partial_\theta \log
    g(\theta, x_t) \right] g(\theta, x_t) \dd x_t \\
  &= \pexp \left[
    f(X_t(\theta)) \partial_\theta \log g(\theta, X_t) \right],
\end{align*}
is obtained by bringing the derivative inside the integral and
multiply and dividing by $g$.  However, this formulation requires
knowledge of the unknown density $g$ and, in practice, pathwise
information is substituted:
$\mathcal{S}_{LR} (t,\theta; f) \approx \pexp[ f(X_t(\theta))
G(\{X_s\}_{0\leq s \leq t})]$.  For both estimators, the application
of these methods to key parameters of interest in the extended
variable formulation of the GLE leads to perturbations in path-space
that are not absolutely continuous, that is, the typical Girsanov
change of measure does not hold in path-space.

An approach that circumvents this lack of a Girsanov transform for
certain parameters, using the Malliavin derivative\cite{Nualart:2006},
first appeared in finance applications for calculating sensitivities,
known as Greeks, related to the pricing of certain securities.
\cite{FournieEtAl:1999oi,FournieEtAl:2001ii,BouchardEtAl:2004ma}
Applied directly to continuous time dynamics, the Malliavin approach
produces unbiased estimators
$\mathcal{S}_M = \pexp [ f(X_T) h(\{X_s\}_{0\leq s\leq T})]$ where $h$
is a non-unique weight that involves a system of auxiliary processes
obtained through Malliavin calculus but that does not depend on
$g$. In particular, for overdamped Langevin dynamics Malliavin weights
are given for sensitivities with respect to parametric forces, that
is, for parameters appearing in the drift term
only.\cite{WarrenAllen:2012mw}
While in principle the Malliavin method applies to other perturbations
that cannot be handled by pathwise and likelihood ratio methods, it
requires a number of auxiliary processes that may scale poorly with
the system size and is not clearly computationally practical for the
extended variable GLE.

We mention that, finite differences using common random numbers have
been employed, based on empirical evidence, for SA with respect to
parametric forces for the LE.\cite{AssarafEtAl:2015sa,
  CiccottiJacucci:1975dc, StoddardFord:1973ne} The sensitivity for
invariant measures for parametric diffusions,
$\partial_\varepsilon (\int_{\mathbf{R}^d} f \dd \pi_\varepsilon)$,
has been considered,\cite{AssarafEtAl:2015sa} and a mathematically
rigorous justification of such objects has been given by other authors
in relation to linear response theory.\cite{HairerMajda:2010rt}
Coupled finite difference methods, similar to the approach outlined
here, have also been applied with success to discrete state space
Markov models, in discrete and continuous time, in chemical kinetics
(chemical reaction networks). For chemical kinetics several couplings
have been demonstrated to reduce the variance of the estimator with
respect to independent sampling.\cite{Anderson:2012fd,
  SrivastavaEtAl:2013cp, WolfAnderson:2015hy, RathinamEtAl:2010sa,
  ArampatzisKatsoulakis:2014go, ArampatzisEtAl:2015as,
  SheppardEtAl:2013spsens} Here, in contrast, we develop a general
framework, at the level of the generators of the coupled SDEs, that
allows us to formulate an optimization problem, locally in time, with
minor assumptions to ensure the correct marginal statistics. That is,
we formulate an associated maximization problem (see Optimization
Problem~\ref{opprob:sigma} in \S \ref{sec:optimal-variance-reduction})
and we show that the problem is satisfied by the common random path
coupling for a large subset of solutions. Further intuition is
developed in the examples of the Ornstein--Uhlenbeck (OU) process and
LE dynamics (see Appendix) for which the optimality of the common
random path coupling can be derived directly without invoking a
localization argument.

In relation to SA, we also mention that information theoretic bounds
can be used to screen parametric
sensitivities.\cite{PantazisKatsoulakis:2013re, PantazisEtAl:2013rn,
  TsourtisEtAl:2015, ArampatzisEtAl:2015pw} In particular, information
theoretic bounds involving the relative
entropy\cite{DupuisEtAl:2015ps} have been used to analyze the
sensitivity of some parameters of interest in Langevin dynamics in
equilibrium and non-equilibrium regimes.\cite{TsourtisEtAl:2015,
  ArampatzisEtAl:2015re} These information theoretic methods are not
goal oriented, that is, the dependence on the observable $f$ is not
explicit. Further they cannot be applied to key parameters of interest
in the extended variable GLE as relative entropy calculations also
require the absolute continuity of the measures arising from the
nominal and perturbed dynamics.

In addition to local SA, the optimally coupled differences are a
useful computational tool for \emph{global SA} and for investigating
\emph{non-local perturbations} in parameters.  In global SA,
elementary effects are used to screen for sensitive
parameters.\cite{Morris:1991fs,CampolongoEtAl:2007sa,SaltelliEtAl:2008gs}
Calculating elementary effects involves sampling a number of finite
difference estimators with various biases and stencils to survey the
space of input parameters. The coupled finite differences might be
used to efficiently build such global sensitivity indices. For
exploring non-local perturbations, a key observation is that the
finite difference estimators proposed are formed by coupling the
nominal and perturbed dynamics and there is no requirement that the
perturbations be local or that the corresponding measures be
absolutely continuous. In \S \ref{sec:sa-number-prony}, we demonstrate
the optimally coupled difference might be used to efficiently analyze,
with respect to independent sampling, the effect of increasing the
\emph{number} of Prony modes active in an extended variable
formulation of GLE dynamics.

The rest of this paper is organized as follows. To set the stage for
our variance reduction technique, we review the errors committed in
estimators for sensitivity indices for static distributions in the
next section. Then we introduce a very general coupling framework and
derive a maximization problem for the variance reduction. In \S
\ref{sec:systems-with-memory} we recall facts about the GLE and
illustrate how the theory presented in \S \ref{sec:effic-finite-diff}
applies to the extended variable formulation, obtaining the optimality
of the common random path coupling for a large subset of solutions. In
\S \ref{sec:numer-exper} we provide numerical experiments involving SA
for GLE that include (1) the sensitivity with respect to the
coefficients of the Prony series approximation, for both convex and
nonconvex potentials, and (2) the sensitivity with respect to the
number of Prony modes, the latter not being formally a sensitivity
index. Finally, in the Appendix, we provide supplemental examples that
help build an intuition for the behavior of coupled finite difference
estimators for other models of interest in the study of particle
dynamics, namely OU processes and the LE.

\section{Efficient finite difference estimators}%
\label{sec:effic-finite-diff}%

In forming the Monte Carlo finite difference estimator for the
sensitivity, the discretization of the derivative results in
systematic error, or bias, while replacing the expected value with a
sample average results in statistical error. We denote the sample
average of $f$, for a sample of size $M$, by
$\hat{f}(X_t) = M^{-1} \sum_{i=1}^M f(X_{i,t})$, where the $X_{i,t}$
are independent for each $i \in \{1, \dots, M\}$.  A measure of the
statistical error committed in computing $\mathcal{S}_\varepsilon$ is
the variance, or more precisely, the standard deviation of the sample
means which is proportional to the square root of the variance.

\subsection{Errors}%
\label{sec:errors}%

To illustrate how these two errors behave, consider for simplicity the
observable that depends on the process at the final time, and define
$\hat{\phi}(\theta) = M^{-1} \sum_{i=1}^M X_{i,T}(\theta)$, a random
variable dependent on the parameter $\theta$. The forward difference
estimator for this observable is
\begin{equation*}
  \mathcal{S}_\varepsilon(T,\theta; \phi) \approx \hat{\Delta}(M,\varepsilon) 
  = \left( \hat{\phi}(\theta+\varepsilon) - 
    \hat{\phi}(\theta) \right) / \varepsilon,
\end{equation*}
where we write $\mathcal{S}_\varepsilon = \Delta(M,\varepsilon)$ to
emphasize the dependence on $M$ and $\varepsilon$ and, in the sequel,
$ \Delta_c$ for the central difference estimator. Note that under
these assumptions, the target is a distribution, that is, there are no
dynamics, and in this setting the following analysis, that gives the
bias and variance of the estimator, is
classical.\cite{Glasserman:2003mc} The expected value of the estimator
is
$\pexp[\hat{\Delta}] =
(\varepsilon)^{-1}(\hat{\phi}(\theta+\varepsilon) -
\hat{\phi}(\theta))$ and if $\hat{\phi}$ is (twice) differentiable in
$\theta$, the bias is given by
\begin{equation}
  \label{eq:bias}
  \bias(\hat{\Delta}) 
  = \pexp [\hat{\Delta} - \hat{\phi}^\prime(\theta)] 
  =  \hat{\phi}^{\prime\prime}(\theta)\varepsilon/2 + O(\varepsilon^2),
\end{equation}
where the last equality can be seen by writing out the Taylor
expansion for $\hat{\phi}(\theta+\varepsilon)$. The variance is
\begin{equation*}
  \pvar[\hat{\Delta}] = \varepsilon^{-2}
  \pvar[\hat{\phi}(\theta+\varepsilon) - \hat{\phi}(\theta)].
\end{equation*}
Assuming the pair $(X_{i,T}(\theta+\varepsilon), X_{i,T}(\theta))$ is
independent of other pairs for each $i \leq M$, then we have that
\begin{equation*}
  \pvar[\hat{\phi}(\theta+\varepsilon) - \hat{\phi}(\theta)] 
  = M^{-1}\pvar[X^1 - X^2]
\end{equation*}
where we define
$(X^1,X^2) = (X_{1,T}(\theta+\varepsilon),X_{1,T}(\theta))$. Thus,
altogether we have
\begin{equation}
  \label{eq:var-diff-est}
  \pvar[\hat{\Delta}] = \varepsilon^{-2}M^{-1} \pvar[X^1 - X^2].
\end{equation}
An analysis of how the variance of this difference depends on
$\varepsilon$ provides insight into a strategy for efficiently
computing the estimator.

From \eqref{eq:var-diff-est}, we see the $\varepsilon$ dependence of
$\pvar[\hat{\Delta}]$ relies upon the $\varepsilon$ dependence of
$\pvar[X^1-X^2]$. If $X^1$ and $X^2$ are independent, then
$\pvar[X^1-X^2] = \pvar[X^1] + \pvar[X^2] \approx 2\pvar[X]$, where
$X$ is related to the distribution of the final time of the nominal
dynamics. This implies $\pvar[X^1-X^2] = O(1)$ and hence
$\pvar[\hat{\Delta}] = O(\varepsilon^{-2}M^{-1})$. In general for
$X^1$ and $X^2$ that are not independent, we have that
\begin{equation*}
  \pvar[\hat{\Delta}] = \varepsilon^{-2}M^{-1} \left( \pvar[X^1] 
    + \pvar[X^2] - 2\pcov[X^1,X^2]\right).
\end{equation*}
Thus, if $X^1$ and $X^2$ are positively correlated, then there is a
net reduction to the variance of the estimator relative to
independently sampling $X^1$ and $X^2$.  For instance, if the
difference $X^1 - X^2$ can be judiciously sampled so that
$\pvar[X^1 - X^2] = O(\varepsilon^2)$, then
$\pvar[\hat{\Delta}] = O(M^{-1})$, asymptotically eliminating the
dependence of the estimator on the bias. For these static
distributions, the well known technique of sampling using common
random numbers (CRN) leads to reduced variance
estimators.\cite{Glasserman:2003mc} We also mention that, for discrete
time Markov chains, the rate of convergence for finite difference
estimators of sensitivities using common random numbers has been
given.\cite{LEcuyerPerron:1994cr,LEcuyer:1992cr} Observe that all of
the error estimates and relations above can be extended from this
simple example with static distributions to the case of dynamics in a
straight forward manner and, in particular, that
\eqref{eq:var-diff-est} remains a quantity of interest for evaluating
the efficiency of the finite difference estimator. Our goal will be to
choose a sampling strategy for \emph{dynamics} that will make the
positive correlation between the distribution of the nominal and
perturbed dynamics at each time step as large as possible while
maintaining the correct marginal statistics.

We remark that, at present we fix a bias and show that the common
random path coupling produces a reduction to the variance
\emph{relative to independent sampling}. The mean squared error (MSE),
formally
\begin{equation*}
  \mse = \pvar + \bias \cdot \bias,
\end{equation*}
represents a balance between the statistical and systematic errors.
While increasing the number of samples $M$ decreases the variance with
no effect on bias, decreasing $\varepsilon$ may increase variance
while decreasing bias. For dynamics, different estimators, for
example, central or forward differences, may have an optimal choice of
bias that balances the two sources of error to achieve a minimal MSE,
as is the case for static distributions.\cite{Glasserman:2003mc}

In the subsequent section, we demonstrate that coupling the nominal
and perturbed dynamics using a common random path is an optimal
strategy for sampling dynamics that reduces the variance of the
estimator $\mathcal{S}_\varepsilon$ relative to independent
sampling. For SA for the extended variable GLE with convex potentials
the common random path coupling leads to substantial reductions,
observed to be $O(\varepsilon^2)$, for sensitivities with respect to
key parameters of interest. In these instances, since the statistical
error scales like the square root of the variance, to reduce the error
by a factor of $10$ for independent sampling with a modest bias of
$\varepsilon = 0.1$ would require adding $M=10^4$ samples, in contrast
to $M=10^2$ samples for the common random path coupling!

\subsection{Coupling dynamics}%
\label{sec:coupling}%

In what follows we provide a very general framework that allows us to
derive a coupling for dynamics that minimizes the variance of the
difference between the nominal and perturbed processes appearing in
equation \eqref{eq:var-diff-est}. We note that this difference need
not be associated with a difference estimator, an aspect that we will
exploit to analyze the sensitivity for non-local perturbations in \S
\ref{sec:sa-number-prony}.

Consider the following pair of SDEs,
\begin{equation}
  \label{eq:sdes}
  \dd X^k_t = b_k(X^k_t) \dd t + \sigma_k(X^k_t) \dd W^k_t, 
\end{equation}
subject to the initial condition $X^k_0 = x^k_0$, for $k=1,2$, where
$X^k_t \in \mathbf{R}$. We assume that for $k \in \{1,2,\}$,
$(W^{k}_t)_{t\geq 0}$ are $\mathbf{R}$-valued independent Wiener
processes, on a given stochastic basis, and that the coefficients
$b_k$ and $\sigma_k$ satisfy the usual properties guaranteeing that
each of the solutions is an It\^{o}
diffusion.\cite{IkedaWatanabe:1989} The infinitesimal generators of
\eqref{eq:sdes} are, respectively,
\begin{equation}
  \label{eq:generator-sdes}
  A_k f(x) = b_k(x) f^\prime(x) 
  + \frac{1}{2} \sigma_k^2(x) f^{\prime\prime}(x), 
\end{equation}
$f \in C^2_0(\mathbf{R})$, where the prime indicates the derivative
with respect to the argument. The generator encodes information about
the statistics of the process.\cite{IkedaWatanabe:1989, Oksendal:2003}

A coupling $Z_t = (X^1_t, X^2_t)$ is produced by considering
\begin{equation}
  \label{eq:coupling}
  \dd Z_t = B(Z_t) \dd t + \Sigma(Z_t) \dd W_t,
\end{equation}
subject to initial conditions $Z_0 = (x^1_0, x^2_0)$, with given
\begin{equation*}
  W_t = \begin{pmatrix} W^1_t \\ W^2_t \end{pmatrix} \quad\text{and}\quad
  B(Z_t) = \begin{pmatrix}b_1(X^1_t)\\ b_2(X^2_t)\end{pmatrix}.
\end{equation*} Here
the diffusion matrix, 
\begin{equation*}
  \Sigma(Z_t) = \begin{pmatrix}
    \Sigma_{11}(Z_t) & \Sigma_{12}(Z_t) \\
    \Sigma_{21}(Z_t) & \Sigma_{22}(Z_t)
  \end{pmatrix},
\end{equation*}
depends on functions $\Sigma_{ij}$, $i,j \in \{1,2\}$, to be
determined. Observe that \eqref{eq:coupling} reduces to
\eqref{eq:sdes} by choosing $\Sigma_{12}= \Sigma_{21} = 0$,
$\Sigma_{11}(Z_t) = \sigma_1(X^1_t)$, and
$\Sigma_{22}(Z_t) = \sigma_2(X^2_t)$. The generator for this extended
system is given by
\begin{equation}
  \label{eq:generator-coup-explicit}
  A f(z) = B(z) \cdot \nabla f(z) 
  + \frac{1}{2} \Sigma(z) \Sigma^\top(z) : \nabla^2 f(z),
\end{equation}
$f \in C^2_0(\mathbf{R}^2)$, where $z = (x_1, x_2)$,
$f: \mathbf{R}^2 \to \mathbf{R}$, and we use the notation
$M : N = \trace (M^\top N)$ for the Frobenius product.

With these ideas in mind, we view $Z_t$ as a coupling of the nominal
and perturbed dynamics in the sensitivity estimator, and, as
foreshadowed in \eqref{eq:var-diff-est} in \S \ref{sec:errors}, seek
to minimize the variance of the difference
\begin{equation}
  \label{eq:diff}
  D(Z_t) = f(X^1_t) - f(X^2_t),
\end{equation}
where $X^1_t$ and $X^2_t$ are solutions of \eqref{eq:sdes} for a given
observable $f$.  In general, this minimization can be achieved locally
in time where the constraints are constructed using
\eqref{eq:generator-sdes} and \eqref{eq:generator-coup-explicit}. For
specific examples (see Appendix) it is possible to obtain the optimal
coupling directly without localizing in time.

A slight modification of the above setting is sufficiently general to
consider the LE and the extended variable GLE, both models that we
consider in the sequel. These two models can be cast as systems of
It\^o diffusions where some components might degenerate in that the
noise term may vanish. Then instead of the pair \eqref{eq:sdes}, which
we view as representing the nominal and perturbed dynamics, we
consider a larger system that decomposes into a system of the nominal
dynamics and of the perturbed dynamics, where some equations are
degenerate diffusions. These ideas will be explored in more detail in
\S \ref{sec:systems-with-memory} after we derive a general formulation
for the optimal coupling for \eqref{eq:sdes}.

\subsection{Optimal variance reduction}%
\label{sec:optimal-variance-reduction}%

To obtain an optimal reduction to the variance of \eqref{eq:diff}, we
place the following constraints on the generator of the coupled
system, namely,
\begin{equation}
  \label{eq:A-constraints}
  \begin{split}
    & Af(x_1,x_2) = A_1 f_1(x_1), \qquad \text{when } f(x_1,x_2) = f_1(x_1), \\
    & Af(x_1,x_2) = A_2 f_2(x_2), \qquad \text{when } f(x_1,x_2) =
    f_2(x_2).
  \end{split}
\end{equation}
These constraints ensure that the marginal statistics of the coupled
system match the statistics of the appropriate diffusion solving
\eqref{eq:sdes}. In particular, for $g(z) = f(x_1)f(x_2)$ such that
$g \in C_0^2(\mathbf{R}^2)$, after some manipulation, the generator
$A$ can be expressed, in terms of the generators
\eqref{eq:generator-sdes},
\begin{equation}
  \label{eq:generator-coup-part}
  \begin{split}
    A g(z) &= (A_1 f(x_1))f(x_2) + (A_2 f(x_2))f(x_1) \\
    &\qquad+ (\Sigma_{11}\Sigma_{21} +
    \Sigma_{12}\Sigma_{22})(z)f^\prime(x_1)f^\prime(x_2),
  \end{split}
\end{equation}
provided that $\Sigma_{11}^2(z) + \Sigma_{12}^2(z) = \sigma_1^2(x_1)$
and $\Sigma_{21}^2(z) + \Sigma_{22}^2(z) = \sigma_2^2(x_2)$ hold for
$z = (x_1,x_2) \in \mathbf{R}^2$ and that the mixed partials of $g$
are equal, $\Sigma_{11}\Sigma_{21}=\Sigma_{21}\Sigma_{11}$, and
$\Sigma_{12}\Sigma_{22}=\Sigma_{22}\Sigma_{12}$.

Next we observe that the variance of \eqref{eq:diff} is equal to
\begin{align*}
  \pvar[D(Z_t)] &= \pvar[f(X^1_t)] + \pvar[f(X^2_t)] \\
                &\qquad+ 2 \pexp[f(X^1_t)]\pexp[f(X^2_t)] \\
                &\qquad- 2\pexp[f(X^1_t)f(X^2_t)].
\end{align*}
In order to minimize the variance, we must maximize the last term in
the above equation. Locally in time, that is, for small perturbations
$\delta t$, we have that,
\begin{align}
  \pexp[&f(X^1_{\delta t})f(X^2_{\delta t})] 
          = \pexp [ g(X^1_{\delta t}, X^2_{\delta t})]\\ \notag
        &= e^{\delta t A} g(X^1_0, X^2_0) \\ \notag
        &= [I + \delta t A + O(\delta t^2)] g(X^1_0,X^2_0) \\ \notag
        &= f(X^1_0)f(X^2_0) + \delta t (A_1 f(X^1_0)) f(X_0^2) \\ \notag
        &\qquad + \delta t (A_2 f(X^2_0)) f(X^1_0) \\ \notag
        &\qquad + \delta t (\Sigma_{11}\Sigma_{21} 
          + \Sigma_{12}\Sigma_{22})(X^1_0,X^2_0) f^\prime(X^1_0)f^\prime(X^2_0) \\
        &\qquad + O(\delta t^2), \label{eq:localization}
\end{align}
where the last equality follows from
\eqref{eq:generator-coup-part}. Using these facts we now state the
following formal optimization problem.

\begin{opprob}
  \label{opprob:sigma}
  The choice of the diffusion matrix $\Sigma$ in \eqref{eq:coupling}
  that minimizes the variance of \eqref{eq:diff} is given by
  \begin{equation}
    \label{eq:general-op-prob}
    \max_\Sigma \mathcal{F}(\Sigma; f) 
    = \max_\Sigma \{(\Sigma_{11}\Sigma_{21} 
    + \Sigma_{12}\Sigma_{22}) (z) f^\prime(x_1) f^\prime(x_2)\},
  \end{equation}
  for all $z = (x_1,x_2) \in \mathbf{R}^2$, under the constraints
  $\Sigma^\top \Sigma \geq 0$ and
  \begin{equation}
    \label{eq:constraint-sigma}
    \begin{split}
      & \Sigma_{11}^2(z) + \Sigma_{12}^2(z) = \sigma_1^2(x_1),\\
      & \Sigma_{21}^2(z) + \Sigma_{22}^2(z) = \sigma_2^2(x_2).
    \end{split}
  \end{equation}
\end{opprob}

To make Optimization Problem \ref{opprob:sigma} tractable, we consider
a restricted subset of all diffusion matrices $\Sigma$. For all
$\Sigma$ of the form
\begin{equation*}
  \Sigma(x_1,x_2) = \begin{pmatrix} 
    \lambda_1 \sigma_1(x_1) & \lambda_2 \sigma_1(x_2)\\
    \lambda_3 \sigma_2(x_1) & \lambda_4 \sigma_2(x_2) 
  \end{pmatrix},
\end{equation*}
we seek to maximize
\begin{equation*}
  \max_\lambda \{(\lambda_1 \sigma_1 \lambda_3 \sigma_2 + \lambda_2 \sigma_1 \lambda_4 \sigma_2)f^\prime(x_1)f^\prime(x_2)\},
\end{equation*}
over $\lambda = (\lambda_1, \lambda_2, \lambda_3, \lambda_4)$, where
the constraints \eqref{eq:constraint-sigma} reduce to
$\lambda_1^2 + \lambda_2^2 = 1$ and $\lambda_3^2 + \lambda_4^2 = 1$.
A solution to this problem is $\lambda_1 = \lambda_3 = \sin(\eta)$,
$\lambda_2 = \lambda_4 = \cos(\eta)$ for all $\eta \in [0,2\pi]$.
Thus we obtain a \emph{family of couplings}
\begin{align*}
  & \dd X^1_t = b_1(X^1_t) \dd t + \sigma_1(X^1_t) (\sin(\eta)\dd W^1_t + \cos(\eta)\dd W^2_t), \\
  & \dd X^2_t = b_2(X^2_t) \dd t + \sigma_2(X^2_t) (\sin(\eta)\dd W^1_t + \cos(\eta)\dd W^2_t),
\end{align*}
for $\eta \in [0,2\pi]$. This coupling is equivalent to generating
approximations with a common Wiener process $(\tilde{W}_t)_{t\geq 0}$
since
$\tilde{W}_t \stackrel{d}{=} \sin(\eta) W^1_t + \cos(\eta)W^2_t$, that
is, they are equal in distribution. Due to the localization argument
in equation \eqref{eq:localization}, this coupling may be sub-optimal
for observables computed over long time horizons. Indeed, for ergodic
systems, observables of trajectories arising from perturbations in the
force field become essentially uncorrelated as the trajectories depart
exponentially as time increases.\cite{StoddardFord:1973ne} For some
explicit examples (see Appendix), one obtains the optimality of the
common random path coupling without requiring a localization
argument. On the other hand, locally for the OU process, LE, and GLE,
we observe that the reduction to the variance of the estimator for
several parameters of interest is on the order of the bias squared;
clearly this coupling must be optimal for the specific numerical
experiments that follow because anything more would be miraculous---we
would have produced a Monte Carlo estimator that could beat Monte
Carlo.

We remark further that the restricted set of diffusion matrices does
not include perturbations of the following form. Consider
$\dd X_t = \sqrt{\mathsf{T}} \dd W_t$ and
$\dd Y_t^\varepsilon = \sqrt{\mathsf{T}} \dd \tilde{W}_t$ for
independent Wiener processes $(W_t)_{t\geq 0}$ and
$(\tilde{W}_t)_{t\geq 0}$. Indeed, $Y^\epsilon_t$ does not define a
local perturbation with respect to $\mathsf{T}$ in precisely the same
manner as
$\dd X^\varepsilon_t = \sqrt{\mathsf{T}+\varepsilon} \dd W_t$. Such
couplings arise in a different context and are natural when the
driving noise is not Brownian but
Poisson.\cite{BenHammoudaEtAl:2016ml} Nevertheless,
$\pcov[Y^\varepsilon_t, X_t] < \pcov[X^\varepsilon_t, X_t]$ and thus
$\pvar[X^\varepsilon_t - X_t] < \pvar[Y^\varepsilon_t - X_t]$, to the
diffusion that is part of our solution set performs better than the
alternative.

In the next section we introduce the GLE, a prototypical system with
memory, and discuss an extended variable formulation which casts the
problem into a form amenable to the preceding theory. We also
introduce some notation and concepts germane to both examples in \S
\ref{sec:numer-exper}, including the technique used for fitting the
Prony series, the normalized velocity autocorrelation function (VACF),
and the integration scheme used.

\section{Systems with memory}%
\label{sec:systems-with-memory}%

\subsection{Extended variable GLE}
\label{sec:extended-variable-gle}
The GLE is a model of anomalous diffusion and subdiffusion, that is,
diffusion where the relationship between the mean squared displacement
(MSD) of the particle and time is no longer linear, that occur in
complex or viscoelastic media typically found in biological
applications.  The GLE includes a temporally non-local drag force and
a random force term with non-trivial
correlations.\cite{Zwanzig:2001ne} The position,
$X^i_t \in \mathbf{R}^d$, and velocity, $V^i_t \in \mathbf{R}^d$, of
particle $i$ with mass $m_i$ at time $t$ are given by the GLE,
\begin{equation}
  \label{eq:gle}
  \begin{split}
    &\dd X^i_t = V^i_t \dd t,\\
    &m_i \dd V^i_t = - \nabla U (X^i_t) \dd t - \int_0^t\!\!\!
    \kappa(t-s) V^i_t \dd s \dd t + F^i(t) \dd t,
  \end{split}
\end{equation}
subject to initial conditions $X^i_0 = x_0$ and $V^i_0 = v_0$, where
$-\nabla U$ is a conservative force and $F^i$ is a random force.  In
the stochastic integro-differential equation for the velocity, the
memory kernel $\kappa$ characterizes the frictional force and, through
the Fluctuation-Dissipation Theorem,
\begin{equation}
  \label{eq:fdt}
  \pexp \left[ F^i(t+s) F^j(t) \right] 
  = k_B \mathsf{T} \kappa(s) \delta_{ij}, 
  \qquad s \geq 0,
\end{equation}
the random force, where $k_B$ is Boltzmann's constant and $\mathsf{T}$
is the absolute (thermodynamic) temperature. This system is
non-Markovian, that is, it has memory; the friction at time $t$ may
have a dependence on the velocity $V(s)$, for $s < t$.

A general strategy for analyzing \eqref{eq:gle} involves mapping the
non-Markovian system onto a Markovian system with suitably many
additional degrees of freedom.\cite{DidierEtAl:2012sc} An extended
variable formulation can often be obtained through physically
realistic assumptions on $\kappa$ that suggest a particular
representation for the memory kernel. For example, when the memory
kernel is posited to have the form of a power law then a positive
Prony series has been identified as a good representation although
more general formulations exist.\cite{Goychuk:2012ve,
  DidierEtAl:2012sc} In general, observe from \eqref{eq:fdt} that
$\kappa$ is the covariance function for the driving noise. Then a
sufficient condition on $\kappa$ for an extended variable formulation
to hold is when the driving noise has a spectral density
$|p(k)|^{-2}$, where $p(k) = \sum_m^{m_1} c_m (-ik)^m $ is a
polynomial with real coefficients and roots in the upper half
plane.\cite{Rey-Bellet:2006oc} A separate topic, not addressed in this
work, is at what level of generality to represent the kernel or
subsequently how to fit the parameters to experimental data. Indeed,
much work has been done in the harmonic analysis and signal processing
literature on fitting exponential functions to data since de~Prony's
classical work.\cite{BeylkinMonzon:2005ap, BeylkinMonzon:2010ae,
  KunisEtAl:2016mv, PottsTasche:2010pe} The important observation here
is that the mapping onto Markovian dynamics yields a system of
(degenerate) It\^o diffusions with a large number of parameters. This
results in systems for which local and global SA are highly relevant
and for which finite differences estimators are useful for SA for all
parameters of interest.

The issue of which representation to use aside, when the memory kernel
can be represented by a positive Prony series,
\begin{equation*}
  \kappa(t) = \sum_{k=1}^{N_k} \frac{c_k}{\tau_k} e^{-t/\tau_k}, 
  \qquad t\geq 0,
\end{equation*}
then the non-Markovian GLE can be mapped into a higher dimensional
Markovian problem in $d N_k$-extended variables. This extended
variable GLE is given by,
\begin{equation}
  \label{eq:extended-variable-gle}
  \begin{split}
    &m \dd V_t = -\nabla U (X_t)\dd t + \sum_{k=1}^{N_k} S^{k}_t \dd t,\\
    &\dd X_t = V_t \dd t,\\
    &\dd S^{k}_t = -\frac{1}{\tau_{k}} S^{k}_t \dd t -
    \frac{c_k}{\tau_k} V_t \dd t + \frac{1}{\tau_k} \sqrt{2 k_B
      \mathsf{T} c_k} \dd W^{k}_t,
  \end{split}
\end{equation}
subject to $X_0 = x_0$, $V_0 = v_0$, and $S^k_0 = s^k_0$, for
independent Wiener processes $(W^k_t)_{t\geq0}$. Here we omit the
obvious extension to a system of many particles in the interest of
brevity. 
In the absence of a conservative force and for the harmonic potential,
$U(X_t) = \omega^2 X_t^2/2$, analytic expressions can be given for the
statistics of the dynamics and for certain observables of interest
including the $\VACF$ and $\MSD$.\cite{DespositoVinales:2009sb,
  KouSie:2004fg, VinalesDesposito:2006ad} For other potentials,
numerical integrators for this system that are stable for a wide range
of parameter values are available and implemented in the LAMMPS
software package.\cite{BaczewskiBond:2013ni} Moreover, these schemes
exactly conserve the first and second moments of the integrated
velocity distribution in certain limits and stably approach the LE in
the limit of small $\tau_k$, the latter of which is a property proven
to hold for the extended variable GLE by other
authors.\cite{OttobrePavliotis:2011aa}

As formulated, \eqref{eq:extended-variable-gle} can be viewed as a
system of (degenerate) It\^o diffusions. Thus, we can form a system of
nominal and perturbed dynamics in the spirit of \eqref{eq:sdes}, for
$k \geq 2$. In addition to any parameters appearing in the potential
and $\mathsf{T}$, we are interested in analyzing the sensitivity with
respect to of $\tau_k$ and $c_k$, $k \in \{1, \dots, N_k\}$. The
pathwise and likelihood ratio methods outlined in the introduction are
not applicable to these latter parameters of interest.  Since in
general the $c_k$ and $\tau_k$ are obtained from experimentally
observed data, it is desirable to analyze the sensitivity of the model
with respect to uncertainties arising from the fitting procedure, for
example, due to errors in measurement or lack of data.

\subsection{Optimal Coupling for extended variable GLE}

Presently we apply the most basic aspects of the theory presented in
\S \ref{sec:effic-finite-diff} to the simple example of an extended
variable GLE with one extended degree of freedom, i.e. one Prony mode,
where the dynamics is subject to a harmonic confining potential with
frequency parameter $\omega$. That is, we consider the system
$\dd Z_t = B Z_t \dd t + \Sigma \dd W_t$ for the coupling
$Z = (X,V,S,\tilde{X},\tilde{V},\tilde{S})$ where $B$ and
$(\Sigma)_{ij} = \sigma_{ij}$ are $6 \times 6$ coefficient matrices to
be determined. Here and below we suppress the extended variable index
and denote the perturbed system variables with tildes for ease of
notation.

An optimal coupling is found by matching the statistics of the
marginals of the coupled process to the statistics of the nominal and
perturbed processes. By writing out the infinitesimal generators of
the corresponding SDEs, this requirement immediately characterizes $B$
and implies that the only nonzero elements of $\Sigma$ are
$\sigma_{33}$, $\sigma_{63}$, $\sigma_{36}$, and
$\sigma_{66}$. Formally, the optimization problem can be stated as
follows.

\begin{opprob}[1-mode GLE with harmonic potential]
  \label{opprob:gle}
  The choice of diffusion matrix $(\Sigma)_{ij} = \sigma_{ij}$ that
  minimizes the variance of $D(Z_t)$ is given by
  \begin{equation*}
    \max_\Sigma \mathcal{F}(\Sigma; f) = \max_\sigma (\sigma_{33}\sigma_{63} + \sigma_{36}\sigma_{66}) \frac{\partial}{\partial x_3 \partial x_6} f(z),
  \end{equation*}
  for $\sigma =(\sigma_{11}, \dots, \sigma_{66})$ for all
  $z \in \mathbf{R}^4$, under the constraints
  $\Sigma^\top \Sigma \geq 0$ and
  \begin{equation*}
    \begin{split}
      & \sigma_{33}^2 + \sigma_{36}^2 = \gamma \sqrt{c}/\tau,\\
      & \sigma_{63}^2 + \sigma_{66}^2 = \gamma
      \sqrt{\tilde{c}}/\tilde{\tau},
    \end{split}
  \end{equation*}
  where $\gamma = \sqrt{2 k_b \emph{\textsf{T}}}$.
\end{opprob}

Thus, for this problem, an optimal family of couplings $Z(\eta)$,
indexed by $\eta \in [0, 2\pi]$, is given by
\begin{equation*}
  \label{eq:gle1-B}
  B =
  \begin{pmatrix}
    L & \boldsymbol{0}\\
    \boldsymbol{0} & \tilde{L}
  \end{pmatrix},
\end{equation*}
with
\begin{equation*}
  \label{eq:gle1-B-L}
  L =
  \begin{pmatrix}
    0 & 1 & 0\\
    -\omega^2 & 0 & 1\\
    0 & -\frac{c}{\tau} & - \frac{1}{\tau}
  \end{pmatrix},
\end{equation*}
and $\Sigma$ with only nonzero elements
$\sigma_{33} = \gamma \sin(\eta) \sqrt{c} /\tau$,
$\sigma_{63} = \gamma \sin(\eta) \sqrt{\tilde{c}} /\tilde{\tau}$,
$\sigma_{36} = \gamma \cos(\eta) \sqrt{c} /\tau$, and
$\sigma_{66} = \gamma \cos(\eta) \sqrt{\tilde{c}} /\tilde{\tau}$,
where $W = (-, -, W^3, -, -, W^6)$ for independent Wiener processes
$(W^3_t)_{t\geq 0}$ and $(W^6_t)_{t\geq 0}$ (here several components
of $W$ are irrelevant due to the zero rows and columns in
$\Sigma$). For each fixed $\eta$, this coupling is equivalent to
choosing a common random path for generating the dynamics of $S$ and
$\tilde{S}$. Extending this optimization problem to an $N_k$-mode GLE
leads to the expected strategy, namely, the common random path
coupling for generating $S^k$ and $\tilde{S}^k$, for each
$k \in \{1, \dots, N_k\}$. Each extended variable requires an
independent common random path for $N_k$ independent Wiener processes
in total, as dictated by \eqref{eq:extended-variable-gle}.

In the remainder of this section we introduce notation and concepts
that are relevant for the numerical experiments in \S
\ref{sec:numer-exper} where we test the variance reduction obtained by
the common random path coupling suggested by the theory above.

\subsection{Fitting Prony series}
\label{sec:fitting-prony}%

\begin{figure}[h]
  \centering
  \includegraphics[width=\shrinkfac\columnwidth]{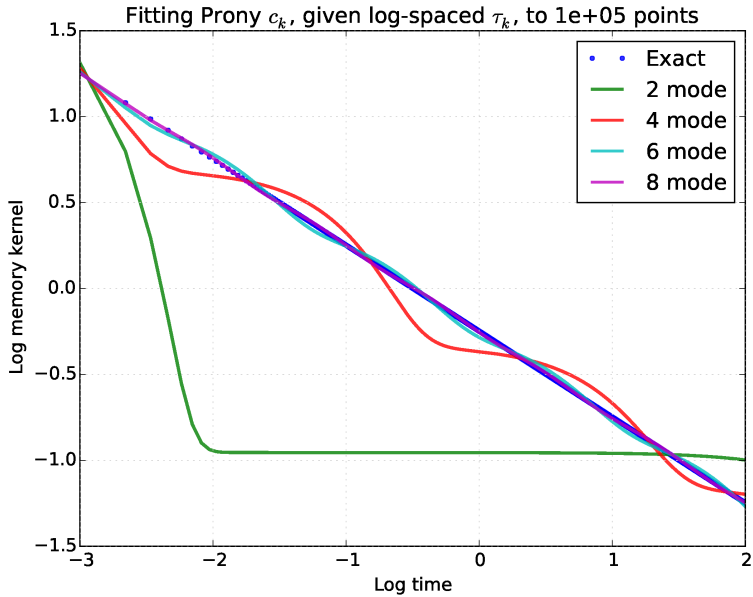}%
  \caption{A least squares fit of the Prony mode coefficients $c_k$,
    $k \in \{1, \dots, N_k\}$, assuming log-spaced $\tau_k$, for each
    of $N_k$ modes. This fit is sufficient to illustrate the variance
    reduction achieved by the common random path coupling.}%
  \label{fig:fit-of-prony-modes}%
\end{figure}

In the numerical experiments that follow, we consider \eqref{eq:gle}
with a power law memory kernel given by
\begin{equation}
  \label{eq:power-law-mk}
  \kappa(t - s) = \frac{\gamma_\lambda}{\Gamma(1-\lambda)} (t - s)^{-\lambda},
\end{equation}
for $\lambda \in (0,1)$ where $\Gamma$ is the gamma function. For
\eqref{eq:power-law-mk}, one can obtain an approximating Prony series
for $N_k$ modes by assuming logarithmically spaced $\tau_k$ and then
fitting the $c_k$ using a least squares method over an interval two
decades longer than the simulation length.\cite{BaczewskiBond:2013ni}
This simplification retains a rich enough family of parameters $c_k$
to illustrate the variance reduction achieved by the common random
path coupling. In Figure~\ref{fig:fit-of-prony-modes}, we illustrate
this fitting procedure for Prony series with $N_k$ modes compared to
measurements of \eqref{eq:power-law-mk} with $\gamma_\lambda = 1.0$
and $\lambda =0.5$.  We choose sufficiently many data points to ensure
a stable least squares approximation.

\subsection{Integration scheme}%
\label{sec:integration-scheme}%

We integrate the system using a modified Verlet method proposed by
other authors, ensuring that the choice of method parameters satisfies
the consistency condition and preserves the Langevin
limit.\cite{BaczewskiBond:2013ni} In many molecular dynamics
simulations, the initial velocity $v_0$, and hence $s_{k,0}$, are
chosen from a thermal distribution. In the numerical experiments
below, the initial conditions for the particle position and velocity
are taken to be definite and $s_{k,0} = 0$ for all $k$. This is done
to minimize the sources of the statistical error thus clarifying the
reporting of deviations in the numerical results and the inclusion of
thermal initial conditions does not pose a challenge to the method.

\subsection{Normalized Autocorrelation Functions}%
\label{sec:vacf}%

\begin{figure}[h]
  \centering
  \includegraphics[height=\shrinkfac\columnwidth, angle=270,
  origin=c]{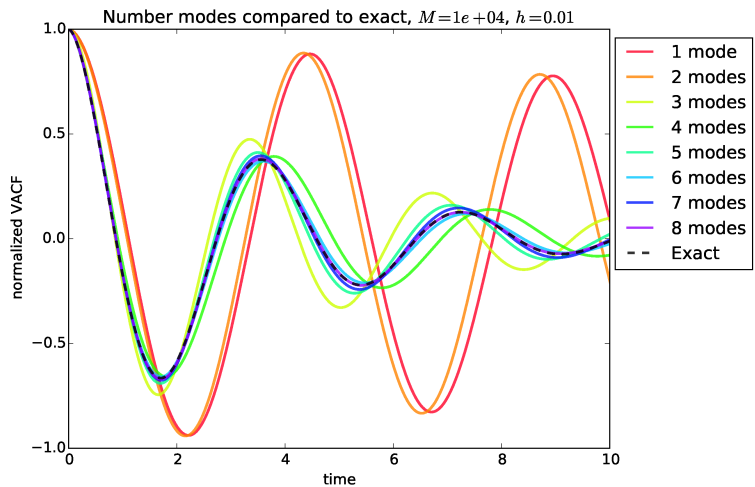}%
  \caption{Small changes to the number of modes leads to qualitatively
    different behavior of the $\overline{\VACF}$ for the GLE with a
    harmonic confining potential.}%
  \label{fig:vacf-number-prony}%
\end{figure}

The results of our numerical experiments are given primarily in term
of normalized autocorrelation functions. Formally, the normalized
$\VACF$ is given by
\begin{equation}
  \label{eq:2}
  \overline{\VACF(t)} = \langle V_t V_0 \rangle/\langle V_0V_0 \rangle 
  = \langle V_t \rangle/v_0,
\end{equation}
where the second equality holds when the initial velocity is
definite. A similar definition is assigned to the normalized position
autocorrelation function ($\PACF$).  For the GLE with a harmonic
confining potential and a power law memory kernel, expressions for the
autocorrelation functions can be given in terms of Mittag-Leffler
functions and their derivatives.\cite{DespositoVinales:2009sb} We
compute the normalized $\VACF$ and $\PACF$ using the integrated
velocity and position distributions and the fast Fourier transform
method.\cite{AllenTildesley:1989cs} Figure~\ref{fig:vacf-number-prony}
illustrates the $\overline{\VACF}$ for models with a varying number of
Prony modes, i.e. extended variables, compared to an asymptotically
exact quantity for the normalized $\VACF$ for the GLE.

\section{Numerical experiments}%
\label{sec:numer-exper}%

The numerical experiments below focus on SA for extended variable GLE
for one particle in one dimension with a power law memory kernel. The
first experiment, in \S \ref{sec:sa-prony-coefficients}, concerns the
sensitivity with respect to the Prony coefficients $c_k$ where the
coefficients are fit using the method described in \S
\ref{sec:fitting-prony}. We observe that the reduction to the variance
of the difference \eqref{eq:diff} for the optimally coupled dynamics
is on the order of the bias squared for convex potentials.

\subsection{Sensitivity with respect to Prony coefficients}%
\label{sec:sa-prony-coefficients}%

\begin{figure}
  \centering%
  \subfloat[$\mathcal{S}_\varepsilon(t,c_1;\overline{\VACF})$,
  $M=10^2$]{%
    \includegraphics[width=\shrinkfactwo\columnwidth]{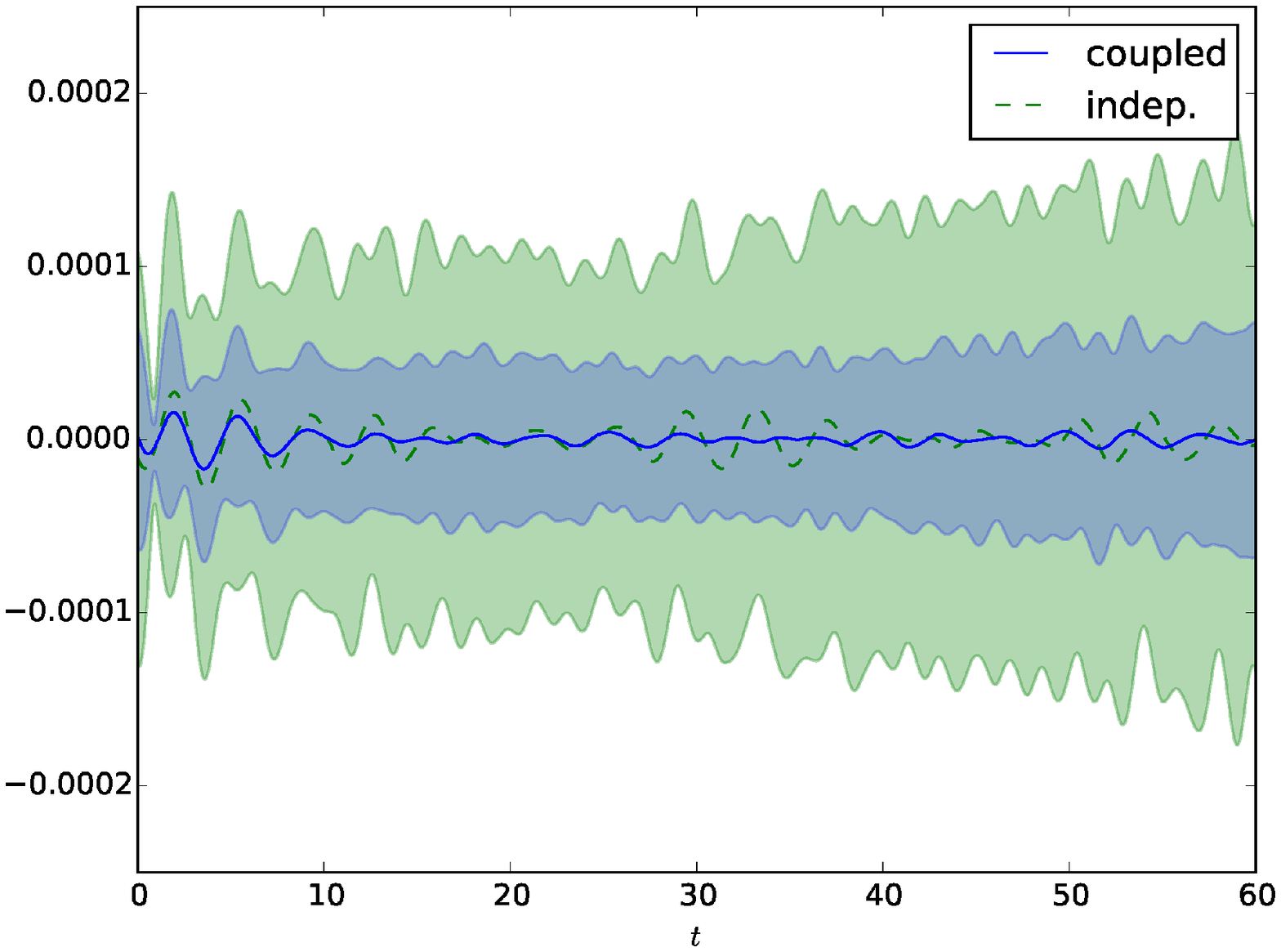}
  }%
  \subfloat[$\mathcal{S}_\varepsilon(t,c_1;\overline{\VACF})$,
  $M=10^4$]{%
    \includegraphics[width=\shrinkfactwo\columnwidth]{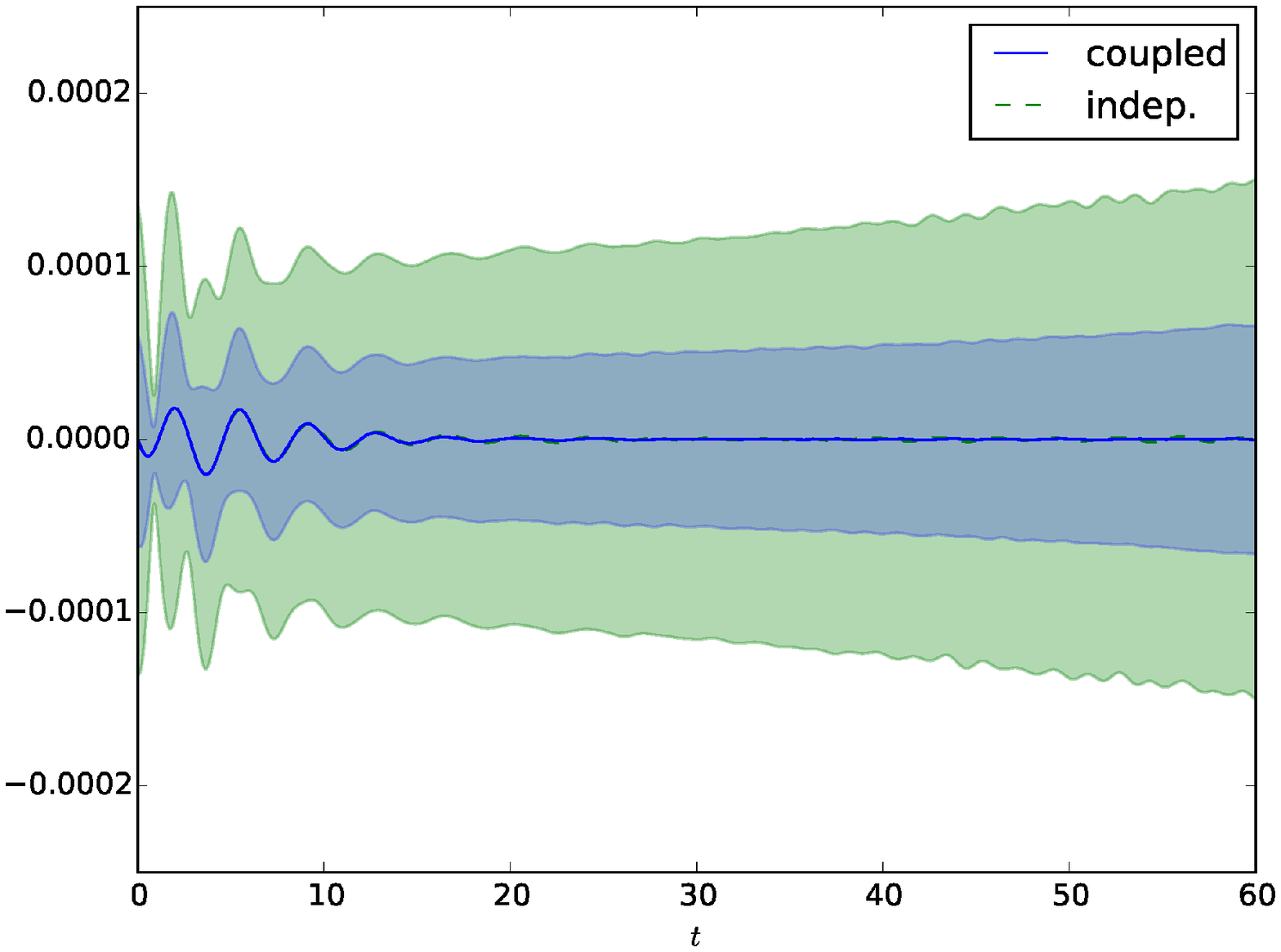}
  }%

  \subfloat[$\mathcal{S}_\varepsilon(t,c_3;\overline{\VACF})$,
  $M=10^2$]{%
    \includegraphics[width=\shrinkfactwo\columnwidth]{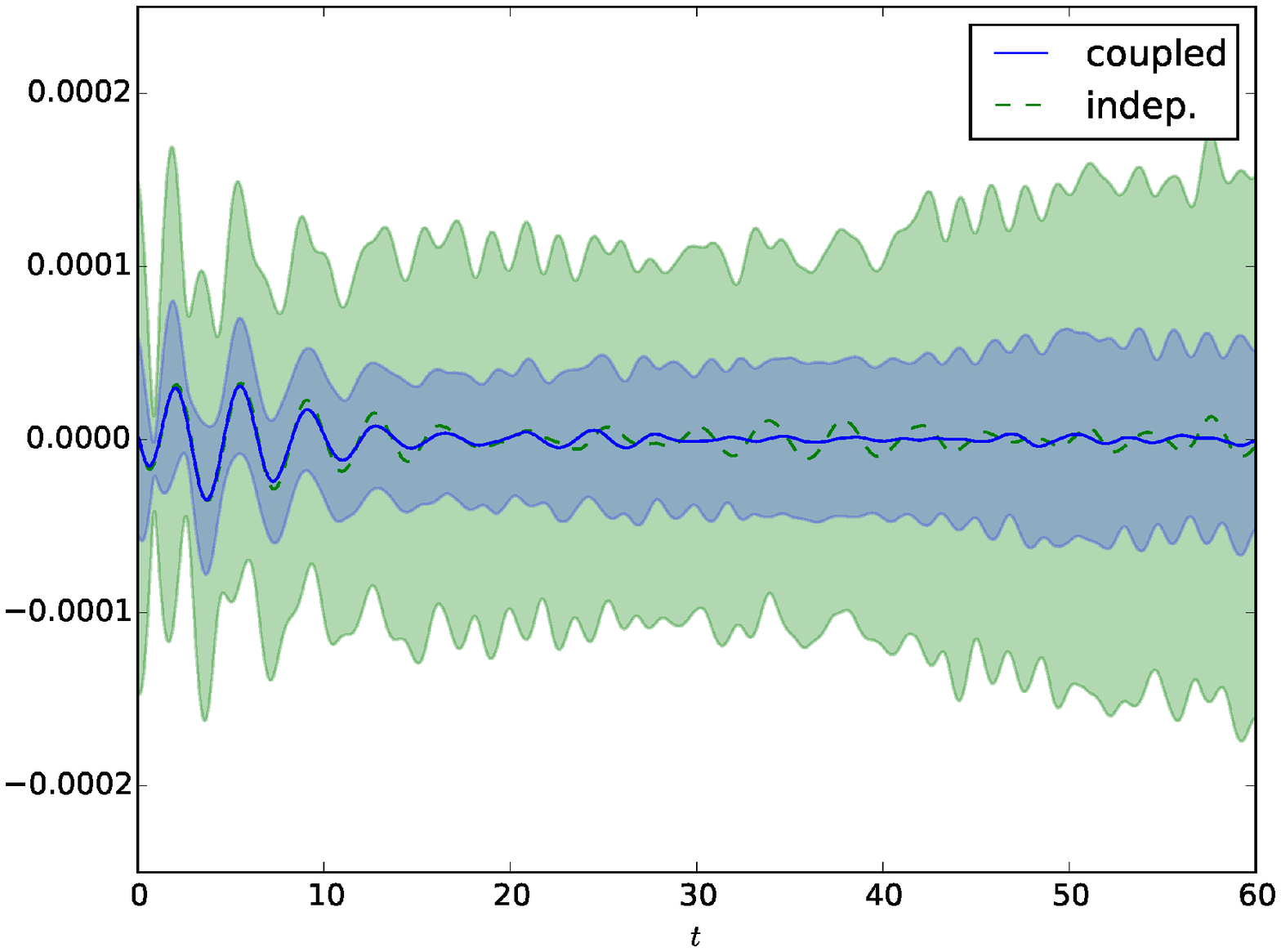}
  }%
  \subfloat[$\mathcal{S}_\varepsilon(t,c_3;\overline{\VACF})$,
  $M=10^4$]{%
    \includegraphics[width=\shrinkfactwo\columnwidth]{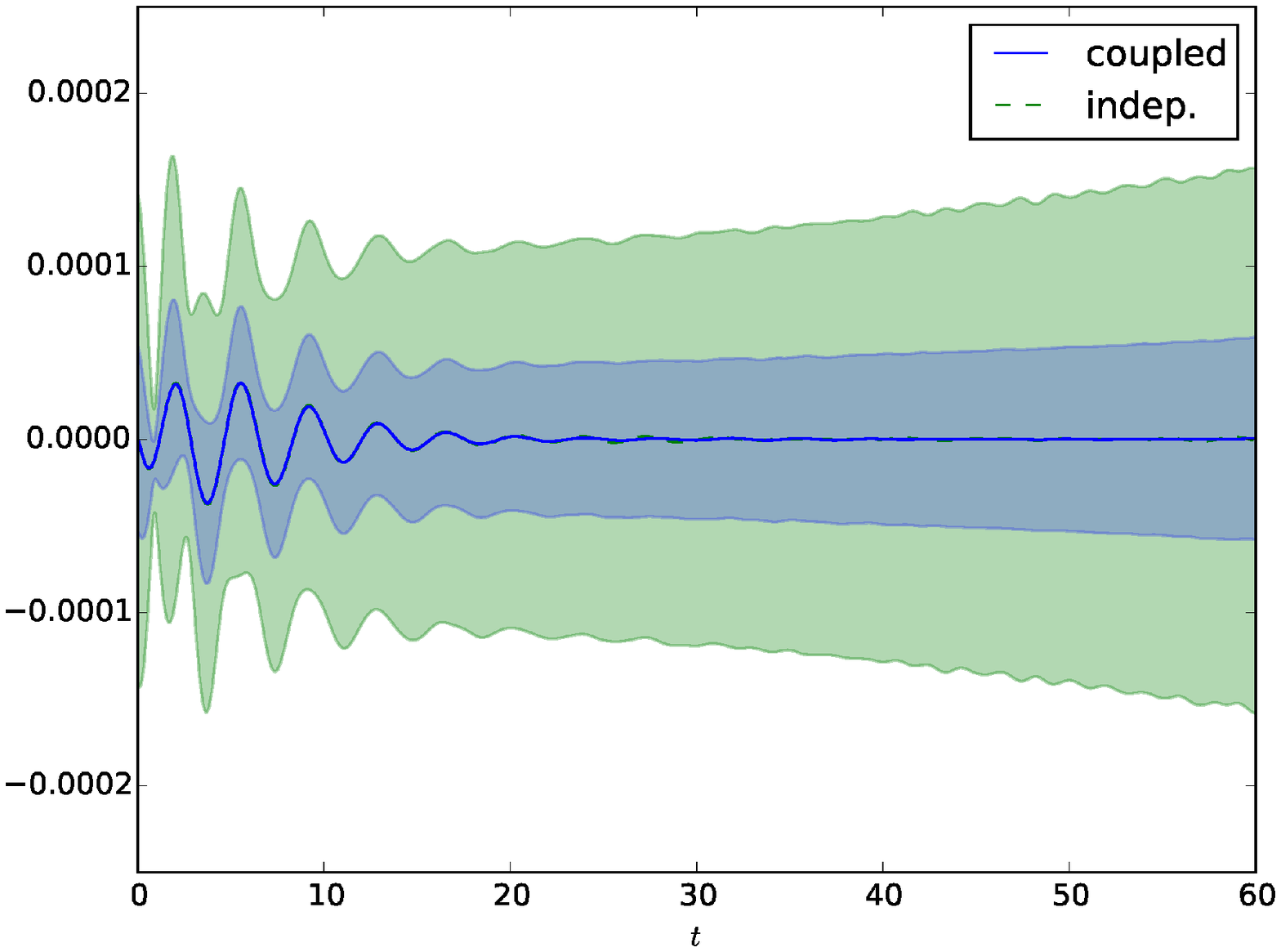}
  }%

  \subfloat[$\mathcal{S}_\varepsilon(t,c_6;\overline{\VACF})$,
  $M=10^2$]{%
    \includegraphics[width=\shrinkfactwo\columnwidth]{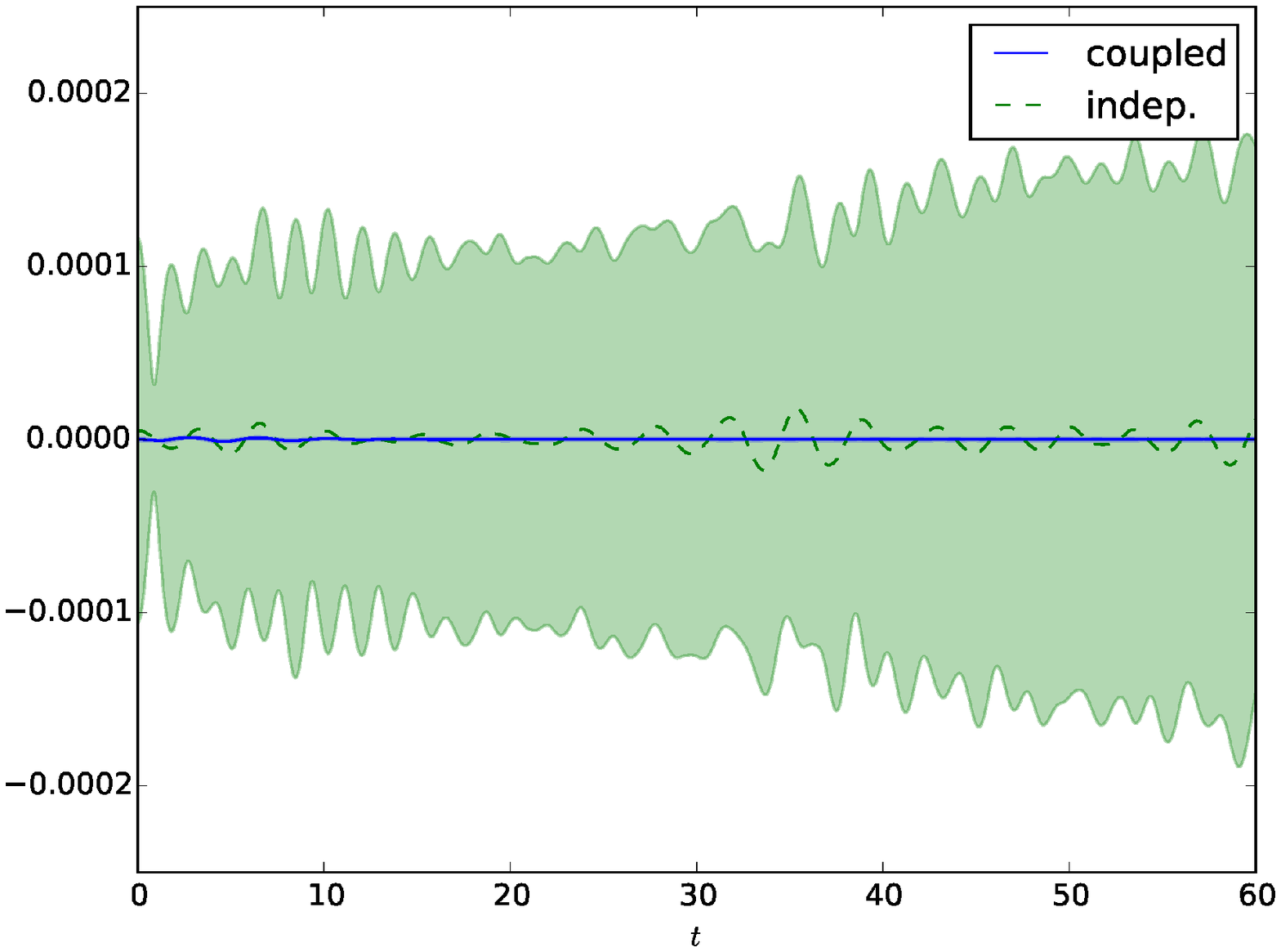}
  }%
  \subfloat[$\mathcal{S}_\varepsilon(t,c_6;\overline{\VACF})$,
  $M=10^4$]{%
    \includegraphics[width=\shrinkfactwo\columnwidth]{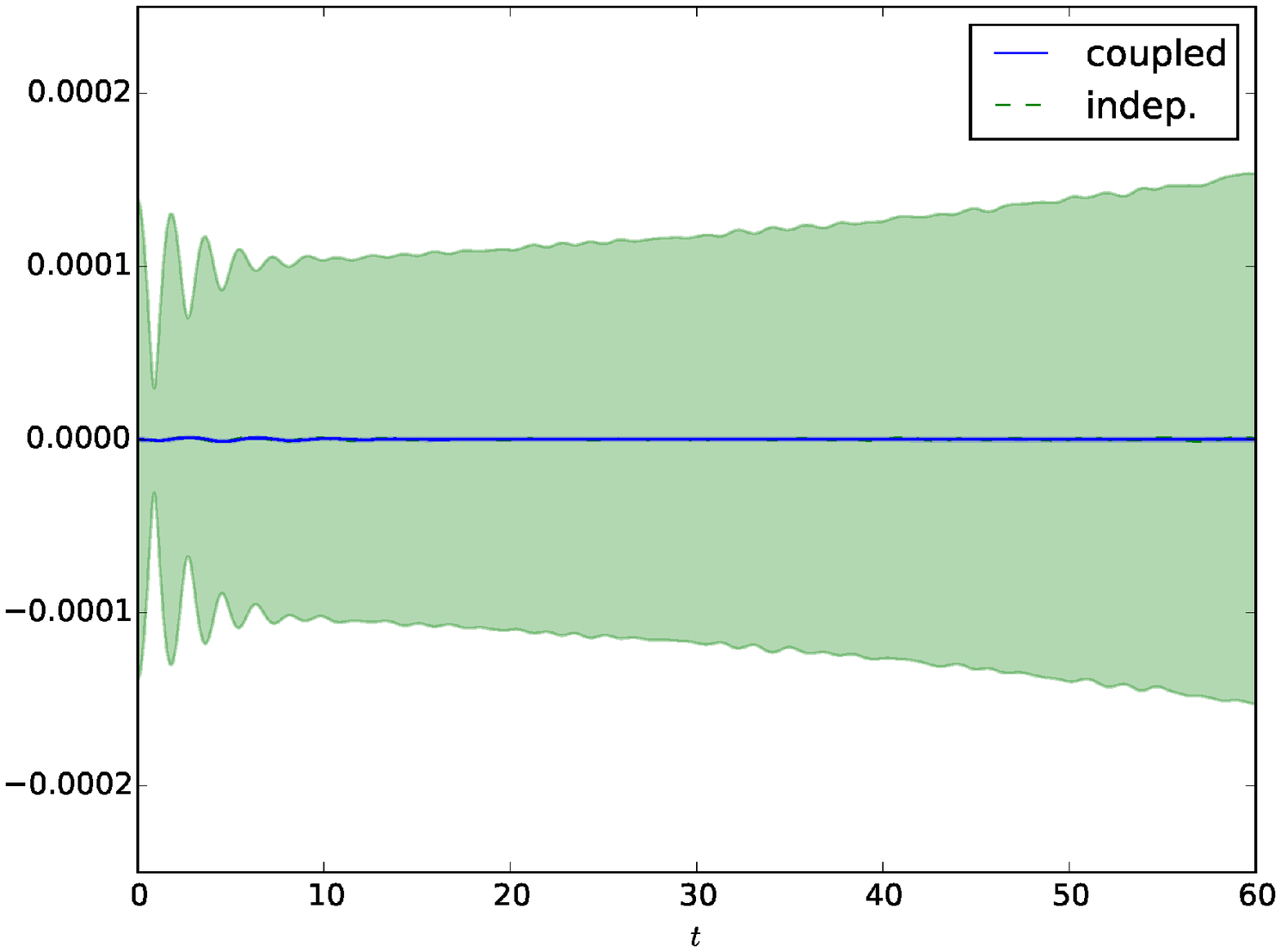}
  }%

  \subfloat[$\mathcal{S}_\varepsilon(t,\omega;\overline{\VACF})$,
  $M=10^2$]{%
    \includegraphics[width=\shrinkfactwo\columnwidth]{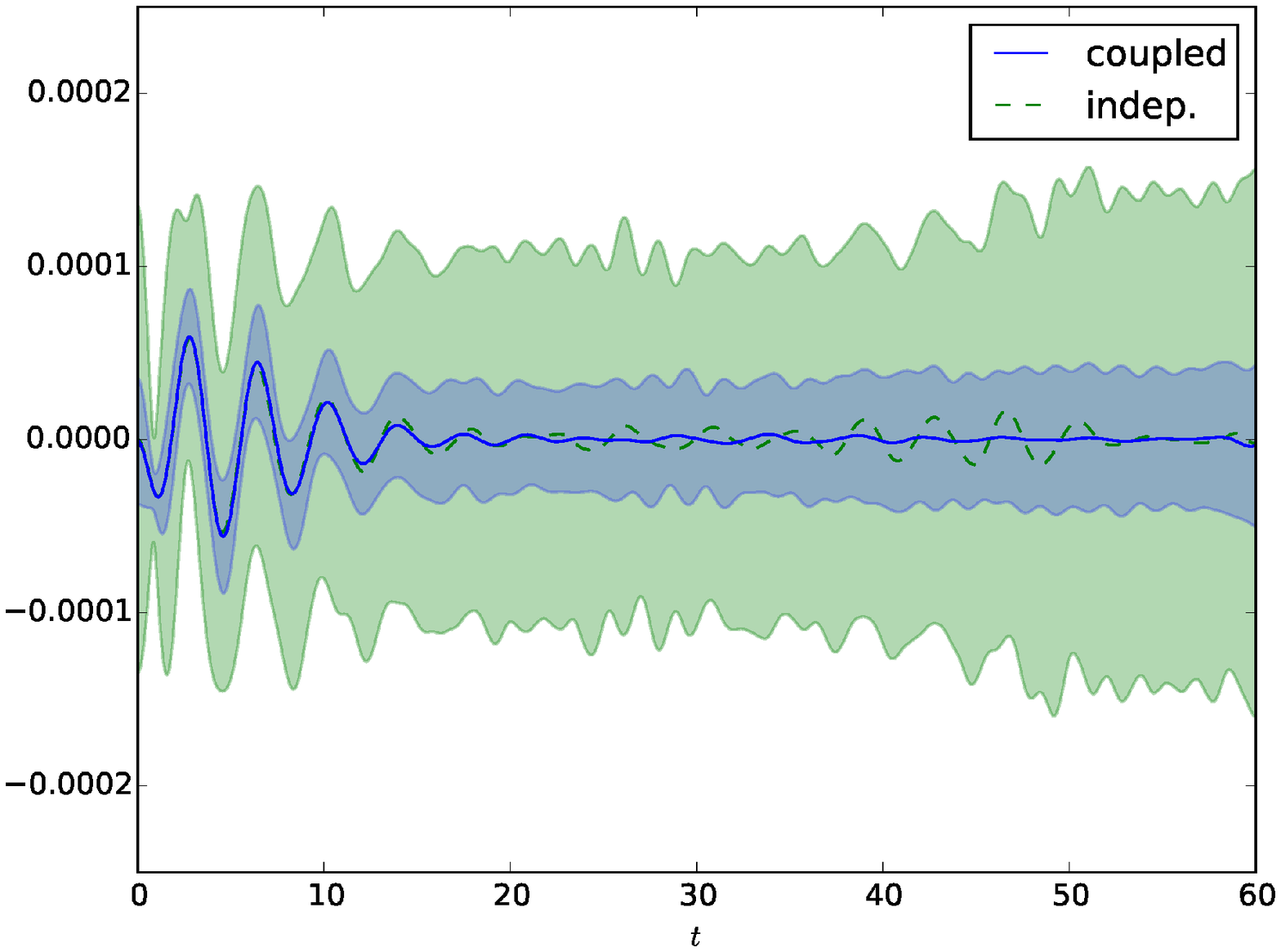}
  }%
  \subfloat[$\mathcal{S}_\varepsilon(t,\omega;\overline{\VACF})$,
  $M=10^4$]{%
    \includegraphics[width=\shrinkfactwo\columnwidth]{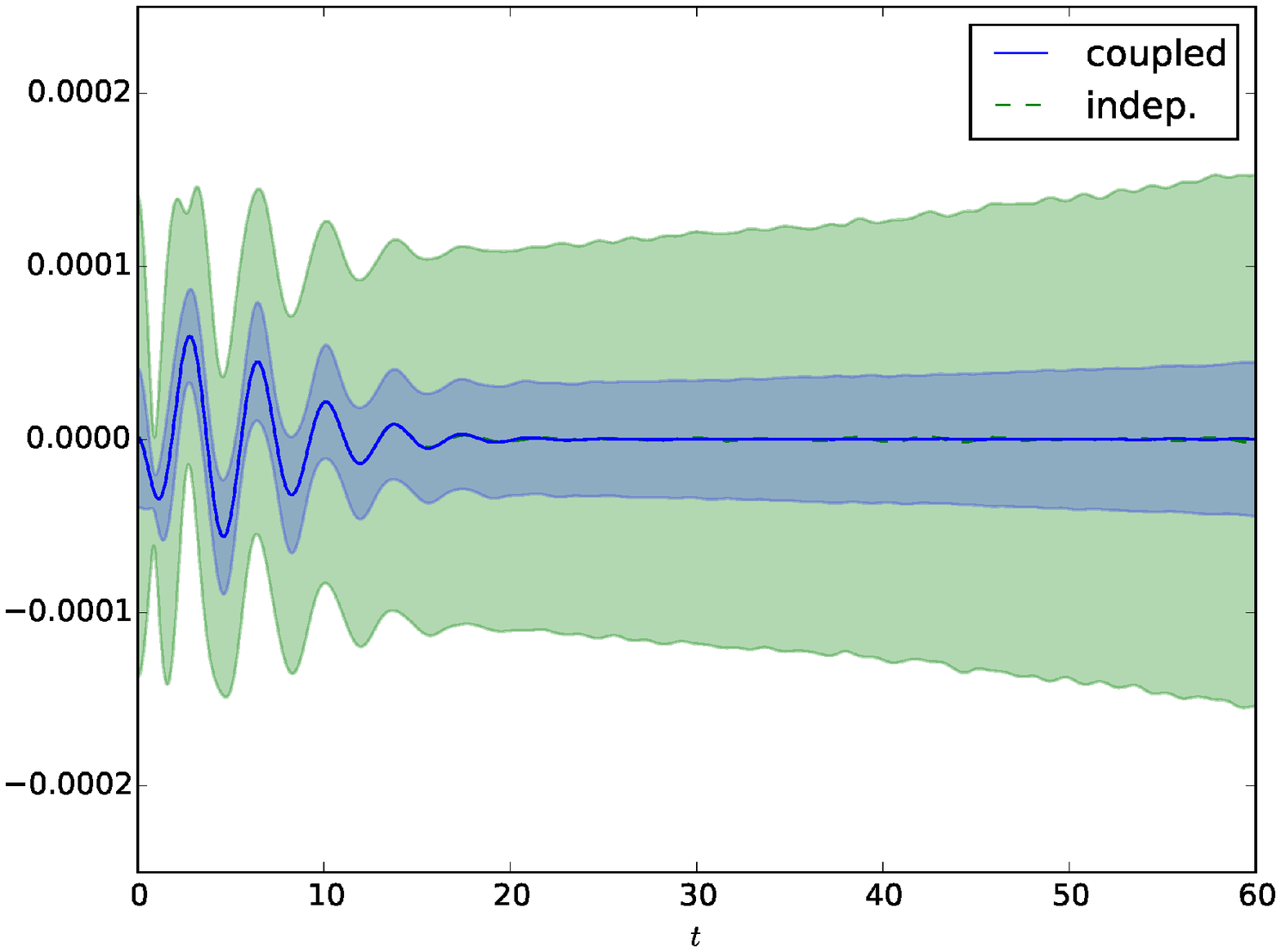}
  }%
  \caption{The computational advantage of the common random path
    coupling is illustrated by the reduced sample size required to
    obtain meaningful estimates for sensitivities. Here we plot the
    sample mean with error bars denoting two standard deviations,
    based on $M=10^2$ samples (left column) and $M=10^4$ samples
    (right column), for various parameters.}%
  \label{fig:sa_stddev_cksr}%
\end{figure}

We begin by computing local sensitivities for the proposed model with
a harmonic confining potential. In particular, we investigate the
sensitivity with respect to the Prony coefficients $c_k$ for
$k \in \{1, \dots, N_k\}$, the harmonic potential frequency $\omega$,
and the temperature $\mathsf{T}$, that is, for a set of parameters
$\theta = (\omega, \mathsf{T}, c_1, \dots, c_{N_k})$. For the
observable $\overline{\VACF(t)} = f(V_t^\theta)$, the Monte Carlo
finite difference estimator
$\mathcal{S}_\varepsilon(t,\theta ; \overline{\VACF}) =
\Delta_c(M,\varepsilon)$ based on the central difference is given by
\begin{equation*}
  \hat{\Delta}_c(M,\varepsilon) 
  = \left( \hat{f}(V_t^{\theta_i+\varepsilon}) - 
    \hat{f}(V_t^{\theta_i-\varepsilon}) \right) / 2\varepsilon,
\end{equation*}
where $V_t^{\theta_i \pm \varepsilon}$ denotes a small $\varepsilon$
perturbation with respect to parameter $\theta_i$ leaving all other
$\theta_j$, $j \neq k$, fixed. We compute $\hat{\Delta}_c$ for a bias
$\varepsilon = 0.01$ for dynamics that are driven by a common random
path and that are driven by independent paths. In
Figure~\ref{fig:sa_stddev_cksr}, we compare the sample mean of
estimators $\mathcal{S}_\varepsilon$, along with one standard
deviation, for various parameters. The key observation here is that
the optimal coupling dramatically reduces the variances of the
difference estimator, relative to the independently sampled dynamics,
even for a modestly sized sample.

\begin{figure}[h]
  \centering
  \includegraphics[width=\shrinkfac\columnwidth]{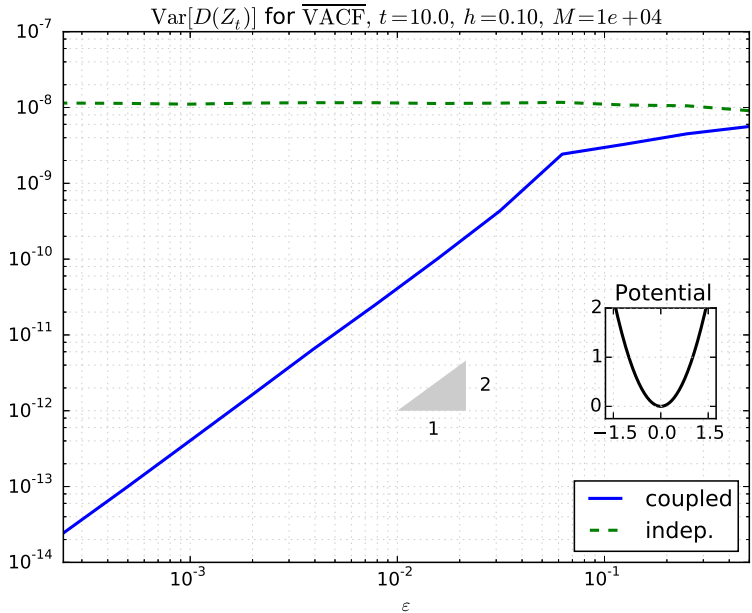}%
  \caption{For $\mathcal{S}_\varepsilon(t=10,c_1;\overline{\VACF})$
    for an $N_k=8$ mode formulation of GLE,
    $\pvar[D(Z_t)] = O(\varepsilon^2)$ for the common random path
    coupling in contrast to $\pvar[D(Z_t)] = O(1)$ for the naively
    sampled independent difference.}%
  \label{fig:sa_c0_varofvacf}%
\end{figure}

The precise nature of the reduction can be deduced by varying
$\varepsilon$ for a fixed index
$\mathcal{S}_\varepsilon(t,\theta_i; \overline{\VACF})$. In
Figure~\ref{fig:sa_c0_varofvacf}, the variance of the difference
\eqref{eq:diff} is compared for dynamics coupled with a common random
path and independent dynamics for
$\mathcal{S}_\varepsilon(t,c_1; \overline{\VACF})$. For the optimally
coupled dynamics, the reduction is $\pvar[D(Z_t)] = O(\varepsilon^2)$,
that is, on the order of the bias squared and, in contrast,
$\pvar[D(Z_t)] = O(1)$ for the difference of the independent
dynamics. Recalling the discussion of errors in \S \ref{sec:errors},
we see that for this example, $\pvar[\hat{\Delta}_c] = O(M^{-1})$ in
the case of the optimally coupled dynamics. That is, the optimal
coupling eliminates the dependence of the variance of the estimator on
the bias, asymptotically, in the case of a convex potential.

\begin{figure}[h]
  \centering
  \includegraphics[width=\shrinkfac\columnwidth]{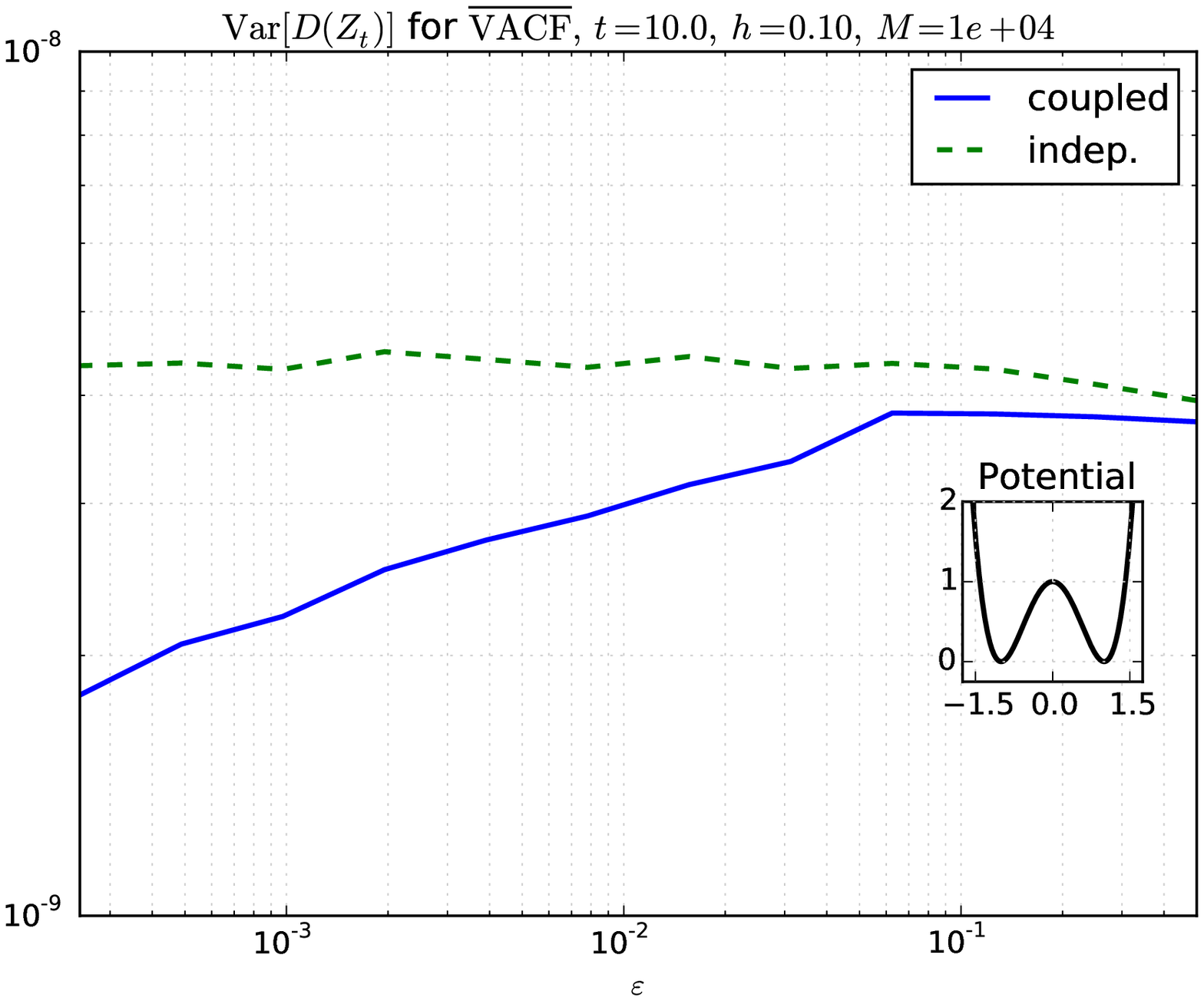}%

  \includegraphics[width=\shrinkfac\columnwidth]{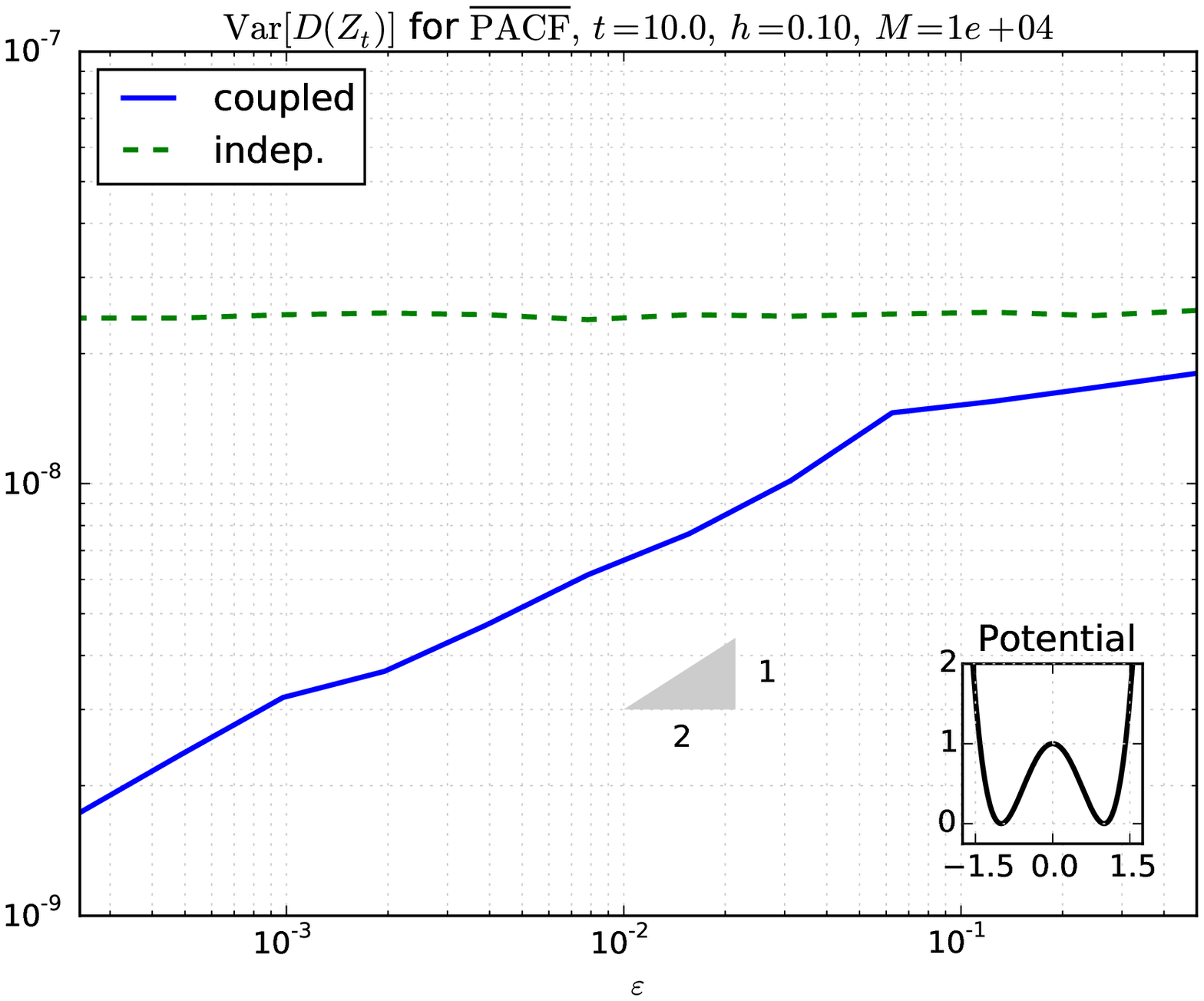}%
  \caption{For simple nonlinear non-convex potentials, there is a net
    reduction to the variance from the common random path
    coupling. Here, for the double well potential,
    $U(X_t) = (1 - X_t^2)^2$, with $k_B\mathsf{T} = 0.5$, the
    reduction is less than $O(\sqrt{\varepsilon})$ for both the
    $\overline{\PACF}$ and $\overline{\VACF}$.}%
  \label{fig:gle_nonconvex}%
\end{figure}

\begin{figure}[h]
  \centering
  \includegraphics[width=\shrinkfac\columnwidth]{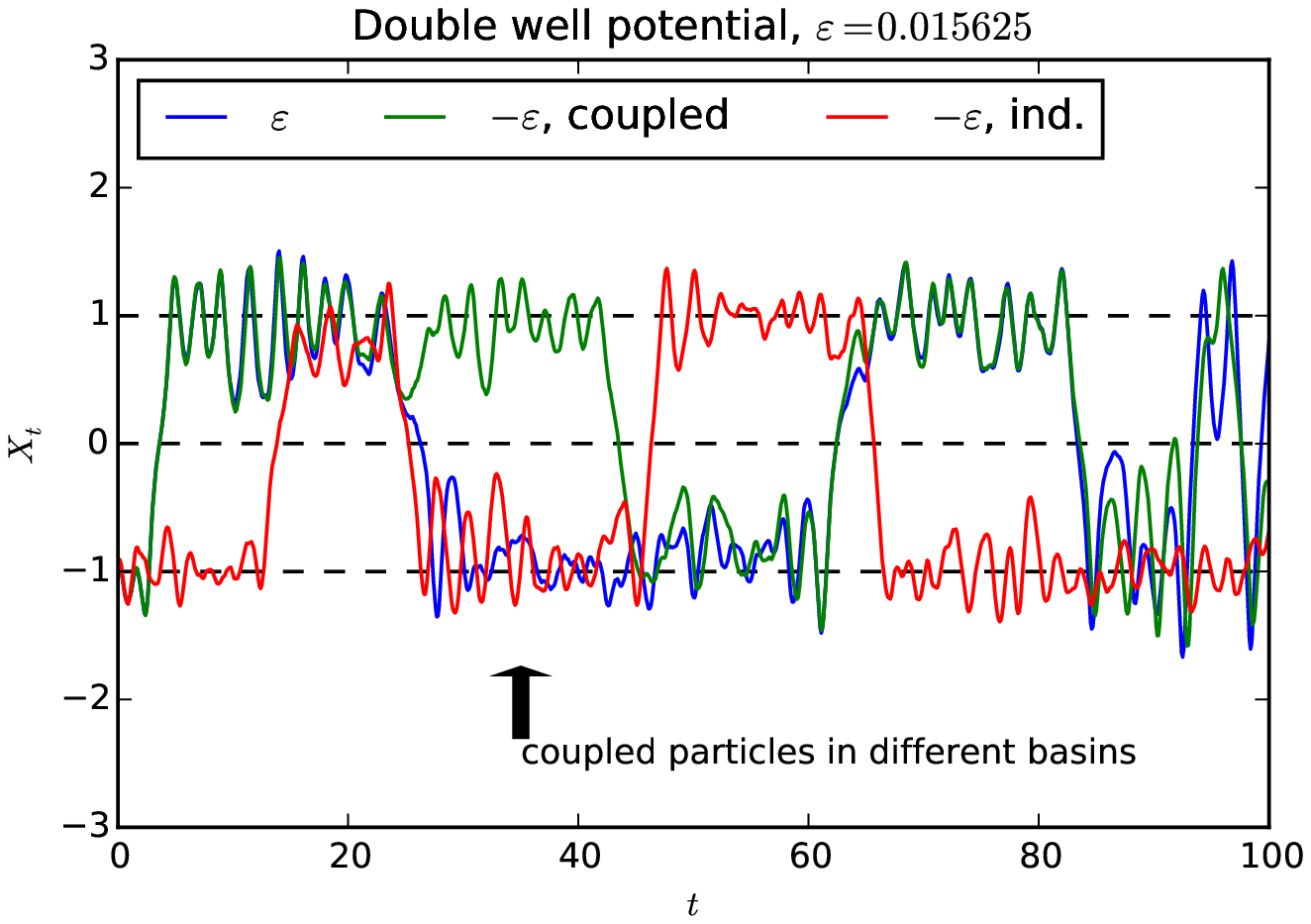}%
  \caption{For nonlinear non-convex potentials, less substantial
    reductions to the variance are observed. For the double well
    potential, the coupled dynamics can be ``pushed'' into distinct
    basis resulting in a higher variance between the coupled paths.}%
  \label{fig:position_data}%
\end{figure}

For nonlinear and non-convex potentials the common random path
coupling reduces the variance of the estimator, although the rate is
not expected to be $O(\varepsilon^2)$. In
Figure~\ref{fig:gle_nonconvex}, the an $N_k=8$ mode formulation of GLE
is considered with a simple double-well potential,
$U(X_t) = (1 - X_t^2)^2$, and $k_B\mathsf{T} = 0.5$ for the
sensitivities $\mathcal{S}_\varepsilon(t=10,c_1; \overline{\VACF})$
and $\mathcal{S}_\varepsilon(t=10,c_1; \overline{\PACF})$. In this
setting, we observe a decay of less than $O(\sqrt{\varepsilon})$ for
both observables. In particular, for the double well potential, the
position time series, see Figure~\ref{fig:position_data}, indicates
that the coupled dynamics can be pushed into distinct basins,
increasing the covariance between the two paths.

For the extended variable GLE with a harmonic potential and a power
law memory kernel, since analytical expressions exists for several
observables of interest including the
$\VACF$,\cite{DespositoVinales:2009sb} the maximum relative error for
approximating the power law memory kernel with a given number of Prony
modes can be computed \emph{a priori}.\cite{BaczewskiBond:2013ni} For
more complicated potentials, exact expressions for observables and
statistics of the dynamics are not available. Further, in reality one
would like to fit the Prony modes to experimentally obtained
data. Such a procedure would likely involve complex inference methods
and a nonlinear fitting to obtain the $\tau_k$ and $c_k$. In such
instances, it would be highly relevant to test the sensitivity of the
fitted parameters.

\subsection{Sensitivity with respect to number of Prony modes}%
\label{sec:sa-number-prony}%

\begin{figure}[h] 
  \centering
  \includegraphics[width=\shrinkfac\columnwidth]{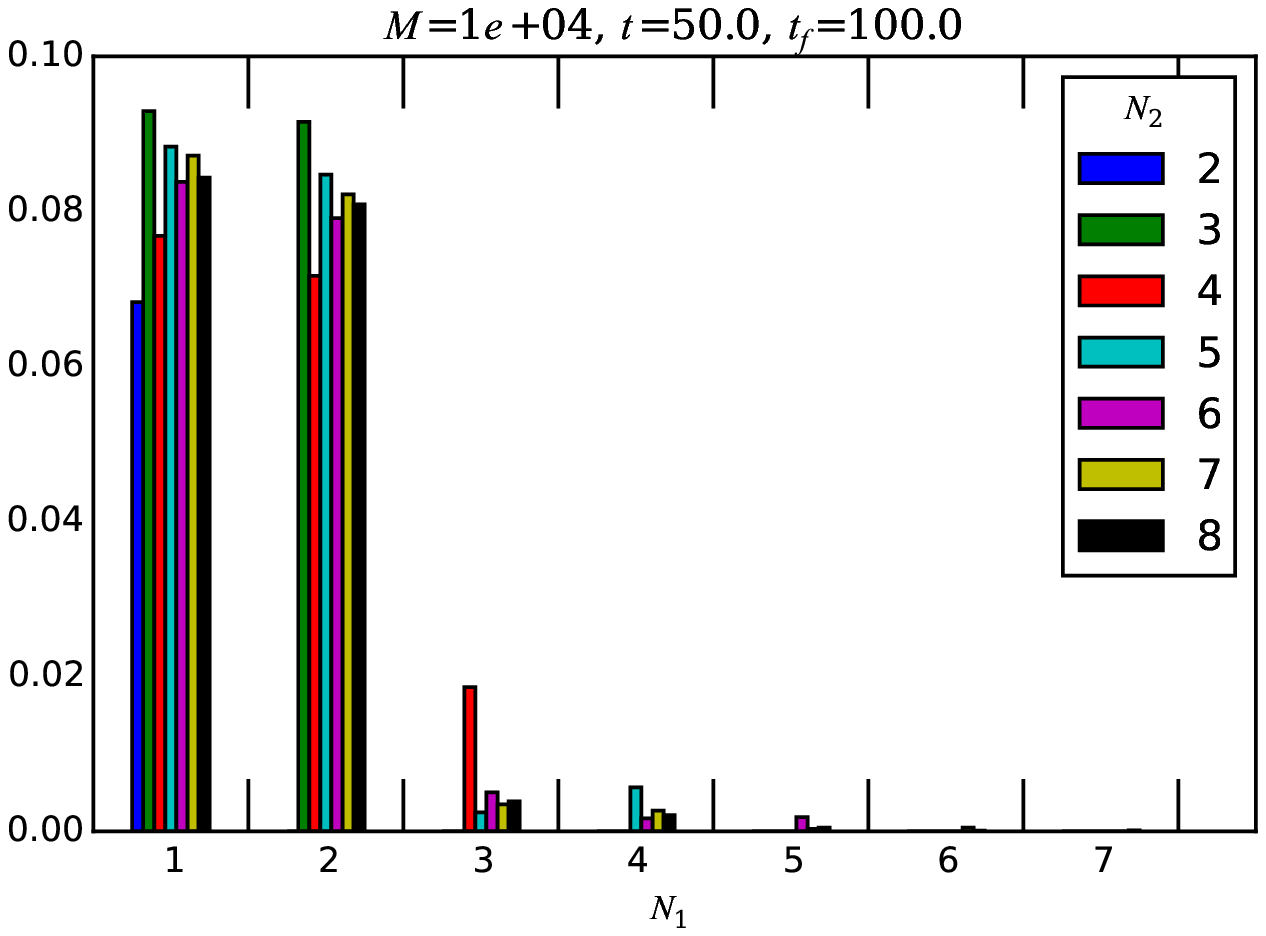}%
  \caption{The non-local sensitivity $\mathcal{S}^*$ gives a
    quantitative characterization of the difference between the
    observed $\overline{\VACF}$ for models with different numbers of
    modes (c.f. Figure \ref{fig:vacf-number-prony}).}%
  \label{fig:si-number-prony}%
\end{figure}

\begin{figure}[h]
  \includegraphics[width=\shrinkfac\columnwidth]{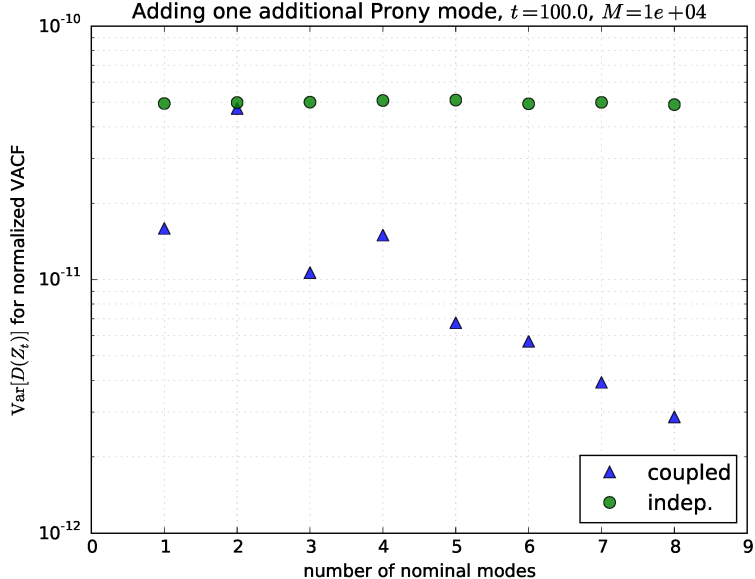}%
  \caption{The common random path coupling is a valid tool for global
    SA, as illustrated by the reduced computational cost in computing
    $\pvar[D(Z_t)]$ where the difference is between a nominal model
    with a fixed number of Prony modes and an alternative model with
    one additional Prony mode. This difference, although it cannot be
    expressed as a derivative, provides a characterization of the
    sensitivity.}%
  \label{fig:sa-vacf-number-prony}%
\end{figure}

From Figure~\ref{fig:vacf-number-prony}, we see that changing the
number of Prony modes has a qualitative impact on the
$\overline{\VACF}$. This motivates the numerical experiment that
follows, where we analyze the effect of increasing the \emph{number}
of Prony modes. That is, for $N_1 < N_2$ consider two systems with
$N_1$ and $N_2$ extended variables, respectively. Given the difference
$D(Z_t) = f(V_t^{N_1}) - f(V_t^{N_2})$, define a sensitivity
\begin{equation*}
  \mathcal{S}^* = \sum_t | \widehat{D(Z_t)} |^2 / \sigma_{Z_t},
\end{equation*}
where the carat denotes a sample mean, $\sigma_{Z_t}$ is the standard
deviation of the associated sample mean $\widehat{D(Z_t)}$, and the
sum is over the space of discrete time points up to a fixed time
$t<T$.  Although this sensitivity is not a sensitivity in the sense of
the gradients introduced previously, $\mathcal{S}^*$ gives a
quantitative characterization of the difference between the two
systems, see Figure \ref{fig:si-number-prony}.

The optimal coupling can be used to reduce the variance of such
non-local perturbations. Here we investigate the difference between a
nominal model with $N_1 = n$, for $n = \{1, \dots, 8\}$, and a
perturbed model with one additional mode $N_2 = N_1+1$. In
Figure~\ref{fig:sa-vacf-number-prony}, we plot the variance of the
difference generated by these nominal and perturbed dynamics for both
the optimally coupled and independent cases, illustrating the reduced
computational cost in sampling the optimally coupled dynamics in
comparison to independent sampling. Here the Prony series are fit
separately for the nominal and perturbed dynamics using the method
outlined in \S \ref{sec:fitting-prony}. Auxiliary variables
$c_{N_1+1} = 0$ and $\tau_{N_1+1} = 1$ are added to the nominal system
so that the vectors for the nominal and perturbed dynamics have the
same size, and then the common random path coupling is naively carried
out for each of the components.

\section{Conclusions}

We develop a general framework for variance reduction via coupling for
goal-oriented SA for continuous time stochastic dynamics. This theory
yields efficient Monte Carlo finite difference estimators for
sensitivity indices that apply to all parameters of interest in an
extended variable formulation of the GLE. Other well known SA
techniques, such as likelihood ratio and pathwise methods are not
applicable to key parameters of interest for this model. These
estimators are obtained by coupling the nominal and perturbed dynamics
appearing in the difference estimator through a common random path and
are thus easy to implement. Strong heuristics are provided to
demonstrate the optimality of the common random path coupling in this
setting. In particular, for the extended variable GLE with convex
potential, the reduction to the variance of the estimator is on the
order of the bias squared, mitigating the effect of the bias error on
the computational cost. Moreover, the common random path coupling is a
valid computational tool in other aspects of UQ including non-local
perturbations and finite difference estimators for global SA.

\begin{acknowledgments}
  The work of all authors was supported by the Office of Advanced
  Scientific Computing Research, U.S.~Department of Energy, under
  Contract No.~DE-SC0010723. This material is based upon work
  supported by the National Science Foundation under Grant
  No.~DMS-1515712.
\end{acknowledgments}

\appendix

\section{Other examples of interest in particle dynamics}%
\label{sec:other-examples}%

\subsection{OU processes}

OU processes are simple and easily analyzed yet are insightful as they
posses several important features: the processes are Markovian,
Gaussian, and stationary under the appropriate choice of initial
conditions. Further, we note that the evolution of the extended
variables in \eqref{eq:extended-variable-gle} is described by an OU
process.

Consider the SDE
\begin{equation}
  \label{eq:ou-sde}
  \dd X_t = \theta (\mu-X_t) \dd t + \sigma \dd W_t,
\end{equation}
subject to the initial condition $X_t = x_0 \sim h$, for a given
distribution $h$, with scalar coefficients $\theta, \sigma >0$, and
$\mu \in \mathbf{R}$. Here $(W_t)_{t\geq 0}$ is a Wiener process on a
given stochastic basis. The solution to \eqref{eq:ou-sde}, given by
\begin{equation}
  \label{eq:ou-sol}
  X_t = x_0e^{-\theta t} + \mu (1-e^{-\theta t}) 
  + \sigma e^{-\theta t} \int_0^t e^{\theta s} \dd W_s,
\end{equation}
for $t \in [0,T]$, is the OU process. This process depends on
parameters $\theta$, $\mu$, $\sigma$, $x_0$, and $h$.

As discussed in \S \ref{sec:effic-finite-diff}, we are interested in
minimizing the variance of \eqref{eq:diff}, where $(X^1_t, X^2_t)$ is
given by the system
\begin{equation}
  \label{eq:ou-sys}
  \begin{split}
    & \dd X^1_t = \theta_1(\mu_1 - X^1_t) \dd t + \sigma_1 \dd W^1_t,
    \qquad X^1_0 = x^1_0 \sim h_1\\
    & \dd X^2_t = \theta_2(\mu_2 - X^2_t) \dd t + \sigma_2 \dd W^2_t,
    \qquad X^2_0 = x^2_0 \sim h_2.
  \end{split}
\end{equation} 
Then
\begin{align*}
  \pvar[D(Z_t)] &= \pvar[f(X^1_t)] + \pvar[f(X^2_t)] \\
                &\qquad + 2 \pexp[f(X^1_t)]\pexp[f(X^2_t)] \\
                &\qquad - 2\pexp[f(X^1_t)f(X^2_t)]\\
                &= \pvar[f(X^1_t)] + \pvar[f(X^2_t)] \\
                &\qquad- 2 \pcov [f(X^1_t), f(X^2_t)],
\end{align*}
and hence to minimize the variance of the difference we seek to
maximize the covariance appearing in the expression above. If $X^1_t$
and $X^2_t$ are independent, that is, they are generated with
independent processes $W^1$ and $W^2$, then the covariance in question
will vanish. If we inject some dependence between $X^1$ and $X^2$ so
that the covariance is nonzero, we find, after cancellation (for
linear $f$), that the covariance is given by
\begin{equation*}
  \pexp \left[ f\left(\sigma_1 e^{-\theta_1 t} \int_0^t e^{\theta_1 s} \dd W_s\right) f\left(\sigma_2 e^{-\theta_2 t} \int_0^t e^{\theta_2 s} \dd W_s\right) \right].
\end{equation*}
This covariance is maximized when the stochastic integral processes
above are dependent, which occurs when the driving processes $W^1$ and
$W^2$ are assumed to be linearly dependent.

We shall look at two concrete observables, to gain intuition on the
variance reduction introduced by the common random path coupling for
the sensitivity with respect to different parameters. For simplicity,
we shall further assume that $x_0^1 = x_0^2$ and that $\mu_1 = \mu_2$
are definite. Then these terms do not play a role since cancellations
occur, for example, when $\pexp[x_0]^2 = \pexp[x_0^2]$.  In these
examples, the coupling with a common random path reduces the variance
in the computation of the central difference estimator by a factor
$O(\varepsilon^2)$ for the sensitivity with respect to $\theta$ and
$\sigma$.

For both observable, and for the sensitivity with respect to $\theta$
and $\sigma$, we find that $\pvar[D(Z_t)] = O(\varepsilon^2)$ when
sampling coupled paths and $\pvar[D(Z_t)] = O(1)$ when sampling
independent paths. Therefore, for standard first order difference
estimators of the sensitivity indices, we have
$\pvar[\hat{\Delta}_c] = O(M^{-1})$, when sampling optimally coupled
paths, but $\pvar[\hat{\Delta}_c] = O(\varepsilon^{-2}M^{-1})$, for
independently sampled paths. For the OU process, the optimal coupling
eliminates the asymptotic dependence of the variance of the estimator
on $\varepsilon$, in contrast to the case of sampling independent
paths.

\subsubsection{Finite time observable}
Consider the finite time observable, $f(X_t) = X_T$ for $T <
\infty$. The expression for the covariance simplifies to
\begin{align*}
  \pcov[&X^1_T, X^2_T] \\
        &= \sigma_1\sigma_2 e^{-(\theta_1+\theta_2)T}\pexp \left[ \int_0^T e^{\theta_1 u} \dd W_u \int_0^T e^{\theta_2 v} \dd W_v \right]\\
        &= \sigma_1\sigma_2 e^{-(\theta_1+\theta_2)T} \int_0^T e^{(\theta_1 + \theta_2)s} \dd  s\\
        &=  \sigma_1\sigma_2 e^{-(\theta_1+\theta_2)T} (e^{(\theta_1 + \theta_2)T} - 1)/(\theta_1 + \theta_2)\\
        &=  \sigma_1\sigma_2 (1 - e^{-(\theta_1 + \theta_2)T})/(\theta_1 + \theta_2).
\end{align*} 
Thus the variance of the difference $D(Z_t)$ converges to a constant,
depending on $\varepsilon$, as $T \to \infty$. As the variance of the
difference does not vanish, the coupling with a common random path is
useful computational technique for all finite times.

Consider now the sensitivity with respect to $\theta$. Then
$\theta_1(\theta,\varepsilon)$ and $\theta_2(\theta,\varepsilon)$ can
be viewed as functions of $\theta$ and $\varepsilon$, i.e.,
$\theta_1 = \theta+\varepsilon$ and $\theta_2=\theta-\varepsilon$ for
the central difference. To determine the asymptotic dependence of
$\pvar[D(Z_t)]$ on $\varepsilon$, we expand the variance of the
difference in a series in $\varepsilon$ at zero. For standard first
order differences (central, forward, and backward), in the case of
independent sampling one finds
\begin{equation*}
  \pvar[X^1_T] + \pvar[X^2_T] = \sigma^2\theta^{-1}
  - \sigma^2 \theta^{-1} e^{-2T\theta} + O(\varepsilon),
\end{equation*}
since $\theta_1(0)=\theta_2(0) = \theta$.  That is, the variance of
the difference is $O(1)$.  In contrast, for sampling with common
random paths, one finds
\begin{equation*}
  \pvar[X^1_T] + \pvar[X^2_T] - 2\pcov[X^1_T,X^2_T] = O(\varepsilon^2).
\end{equation*}
A similar story holds for the sensitivity with respect to $\sigma$:
$\pvar[D(Z_t)] = O(1)$ for independent sampling and
$\pvar[D(Z_t)] = O(\varepsilon^2)$ for sampling with common random
paths, when using standard first order differences.

\subsubsection{Time average observable}
Next we consider the time average observable defined by
$\overline{X} = T^{-1} \int_0^T X_s \dd s$. Once again, we wish
to investigate the dependence of $\pvar[D(Z_t)]$ on $\varepsilon$ for
the case of coupled paths and independent sampling. The expression for
the covariance in this instance is
\begin{widetext}
  \begin{align*}
    \pcov[\overline{X^1}, \overline{X^2} ] 
    &=
      \pexp \left[ T^{-1} \int_0^T \sigma_1 e^{-\theta_1 s} \int_0^s e^{\theta_1 u} \dd W_u \dd s \;T^{-1} \int_0^T \sigma_2 e^{-\theta_2 t} \int_0^t e^{\theta_2 v} \dd W_v \dd t \right]\\
    &= \sigma_1\sigma_2T^{-2} \int_0^T \int_0^T e^{-\theta_1 s - \theta_2 t} \pexp \left[ \int_0^t e^{\theta_1 u} \dd W_u \int_0^s e^{\theta_2 v} \dd W_v \right] \dd s \dd t \\
    &= \sigma_1\sigma_2T^{-2} \int_0^T \int_0^T e^{-\theta_1 s - \theta_2 t} \int_0^{s\wedge t} e^{(\theta_1+\theta_2)r} \dd r \dd s \dd t\\
    &= \sigma_1\sigma_2T^{-1} \left(\left(\theta_1\theta_2 +\theta_2^2\right)^{-1} + \left(\theta_1^2+\theta_1\theta_2\right)^{-1}\right) + O(T^{-2}).
  \end{align*}
\end{widetext}
First we look at the sensitivity with respect to the parameter
$\theta$. As in the case for the finite time observable, we expand
$\pvar[D(Z_t)]$ in a series in $\varepsilon$ at zero. For standard
first order differences this yields
\begin{align*}
  \pvar[D(Z_t)] &= 2\sigma^2 T^{-1} \theta^{-2} - 3\sigma^2T^{-1}\theta^{-3} \\ & + 4\sigma^2 T^{-2} \theta^{-3} e^{-T\theta} - \sigma^2 T^{-2}\theta^{-3} e^{-2T\theta} + O(\varepsilon),
\end{align*}
for independently sampled paths.  Working in a similar fashion, we
find in contrast that,
\begin{equation*}
  \pvar[D(Z_t)] = \varepsilon^2 \left( 4\sigma^2 T^{-1} \theta^{-4} + O(T^{-2})\right) + O(\varepsilon^4),
\end{equation*}
for the coupled paths. For the sensitivity with respect to $\sigma$,
the story is the same. The independently sampled paths behave like
\begin{equation*}
  \pvar[D(Z_t)] = 2\sigma^2 T^{-1}\theta^{-2} + \varepsilon^2 2T^{-1}\theta^{-2} + O(T^{-2}) \left(1 + \varepsilon^2\right)
\end{equation*}
and the coupled paths behave like
\begin{equation*}
  \pvar[D(Z_t)] = \varepsilon^2 \left( c  T^{-1} \theta^{-2} + O(T^{-2}) \right),
\end{equation*} 
for a constant $c$.
\begin{figure}[h]
  \centering
  \includegraphics[width=\shrinkfac\columnwidth]{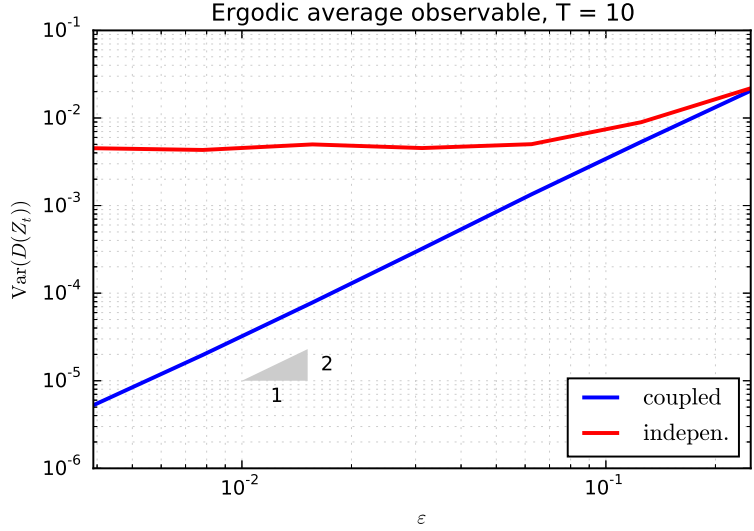}
  \caption{For the OU process, the variance of the estimator for the
    sensitivity with respect to $\sigma$ for the time averaged
    observable is $O(1)$ for the independently sampled difference and
    $O(\varepsilon^2)$ for the optimally coupled difference. Here
    consider an OU process with parameters $\theta = 1$, $\mu = 1.2$,
    $\sigma=0.3$, and $x_0 = 2$ and compute the average up to time
    $T=10$. Each variance is based on $M=10^3$ independent samples of
    optimally coupled $\AR(1)$ processes and independent $\AR(1)$
    processes.}
  \label{fig:coupVarIsEpsSqrd}
\end{figure}
In Figure~\ref{fig:coupVarIsEpsSqrd}, we observe the theoretically
obtained values for the reduction to the variance for the sensitivity
with respect to $\sigma$, of an OU process with parameters
$\theta = 1$, $\mu = 1.2$, $\sigma=0.3$, and $x_0 = 2$. The time
average average is computed up to time $T=10$ and each variance is
computed using $M=10^3$ independent samples of an optimally coupled
$\AR(1)$ processes and an independent $\AR(1)$ processes.

\subsection{Langevin dynamics}

We consider the LE with particle mass $m=1$,
\begin{align*}
  &\dd X_t = V_t \dd t,\\
  &\dd V_t = - \omega^2 X_t \dd t - \beta V_t \dd t + \sqrt{2\beta k_B \mathsf{T}}\dd W^1_t,
\end{align*}
for $t \in [0, T]$, subject to $X_0 = x_0$ and $V_0 = v_0$, where
$(W^1_t)_{t\geq 0}$ is a Wiener process. This system can be written as a two-dimensional OU process
$Y_t = (X_t, V_t)$, given by
\begin{equation}
  \label{eq:langevin}
  \dd Y_t = - B Y_t \dd t +
  \Sigma \dd W_t, \qquad Y_0 = (x_0,v_0),
\end{equation}
for $W_t = (0,W^1_t)$ with coefficient matrices
\begin{equation*}
  B = \begin{pmatrix} 
    0 & -1 \\ 
    \omega^2 & \beta
  \end{pmatrix} \qquad\text{and}\qquad \Sigma = \begin{pmatrix}
    0 & 0 \\
    0 & \sqrt{2\beta k_B \mathsf{T}}
  \end{pmatrix}.
\end{equation*}
The general solution to \eqref{eq:langevin} is given by
\begin{equation}
  \label{eq:langevin-sol}
  Y_t = e^{-Bt} Y_0 + \int_0^t e^{-B(t-s)}\Sigma \dd W_s,
\end{equation}
for $t \in [0,T]$, where, for this example, $e^{-Bt}$ can be written
as (except in the critically damped case) in a closed form in terms of
the eigenvalues of $B$:
$\mu_1 = \beta/2 + \sqrt{\beta^2/4 - \omega^2}$ and
$\mu_2 = \beta/2 - \sqrt{\beta^2/4 - \omega^2}$.\cite{Nelson:1967}
That is, the position and velocity are given component-wise by
\begin{align*}
  X_t &= \mu^{-1} \Big( x_0(\mu_1 e^{-\mu_2 t} - \mu_2 e^{-\mu_1 t}) + v_0 (e^{-\mu_2 t} - e^{-\mu_1 t}) \\
      & + \sqrt{\gamma (\mu_1+\mu_2)}
        \int_0^t (e^{-\mu_2 (t-s)} - e^{-\mu_1 (t-s)})\dd W_s \Big),\\
  V_t &= \mu^{-1} \Big( x_0 \omega^2(e^{-\mu_1 t} - e^{-\mu_2 t}) + v_0 (\mu_1 e^{-\mu_1 t} - \mu_2 e^{-\mu_2 t}) \\
      & + \sqrt{\gamma (\mu_1+\mu_2)} \int_0^t (\mu_1 e^{-\mu_1 (t-s)} - \mu_2 e^{-\mu_2 (t-s)}) \dd W_s \Big),
\end{align*}
for $\mu^{-1} = (\mu_1 - \mu_2)$ and $\gamma = 2k_B \mathsf{T}$.  We
shall further assume, for simplicity, that both $x_0$ and $v_0$ are
definite.

For the Langevin dynamics we form the coupled system
$Z_t = (Y_t, \tilde{Y}_t)$ where $\tilde{Y}$ solves
\eqref{eq:langevin} for $\tilde{B}$ and $\tilde{\Sigma}$ depending
upon perturbed parameters (also denoted with tildes in the sequel) and
with an independent Wiener process $\tilde{W}$. Once again, we are
interested in minimizing the variance of the difference
$D(Z_t) = f(Y_t) - f(\tilde{Y}_t)$, for linear observables $f$.  Note
that $D(Z_t)$ is a vector quantity (i.e., $\pvar[D(Z_t)]$ is the
variance-covariance matrix),
\begin{equation*}
  \pvar[D(Z_t)] = \pvar[f(Y_t)] - 2 \pcov [f(Y_t), f(\tilde{Y}_t)] + \pvar[f(\tilde{Y}_t)] \notag
  \label{eq:langevin-Phi}
\end{equation*}
for $f(Y_t) = (f(X_t), f(V_t))$, where
$\pcov [f(Y_t), f(\tilde{Y}_t)]$ has components
$\pcov[f(X_t),f(\tilde{X}_t)]$, $\pcov[f(V_t),f(\tilde{V}_t)]$, and
cross terms
$\frac{1}{2}(\pcov[f(V_t),f(\tilde{X}_t)] + \pcov[f(X_t),
f(\tilde{V}_t)])$.
This covariance is zero when $Y_t$ and $\tilde{Y}_t$ are independent
and can be maximized when $Y_t$ and $\tilde{Y}_t$ are linearly
dependent, which is equivalent to generating $Y_t$ and $\tilde{Y}_t$
using common random paths $W_t = \tilde{W}_t$. Next we investigate the
asymptotic dependence of $\pvar[D(Z_t)]$ on $\varepsilon$ for two
observables, related to a finite time horizon and a time average, for
sensitivities with respect to $\beta$.

\subsubsection{Finite time observable}
\label{eg:finite-time-observable}

Consider the finite time observable $f(Y_t) := Y_T$. Using the
component wise expression above, the covariance term related to the
positions can be expressed in terms of the eigenvalues of the drift
matrices for the nominal and perturbed systems. That is, we let
$\pcov[X_T, \tilde{X}_T] = \phi (\mu_1, \mu_2, \tilde{\mu}_1,
\tilde{\mu}_2) $ where
\begin{align*}
  \phi (\mu_1, \mu_2, \tilde{\mu}_1, \tilde{\mu}_2) 
  = \frac{\gamma \sqrt{(\mu_1+\mu_2)(\tilde{\mu}_1 + \tilde{\mu}_2)}}{(\mu_1 - \mu_2)(\tilde{\mu}_1 - \tilde{\mu}_2)}
    &\Big(\frac{1 - e^{-(\mu_1+\tilde{\mu}_1)T}}{\mu_1 + \tilde{\mu}_1} - \frac{1 - e^{-(\mu_1+\tilde{\mu}_2)T}}{\mu_1 + \tilde{\mu}_2} \\ & - \frac{1 - e^{-(\mu_2+\tilde{\mu}_1)T}}{\mu_2 + \tilde{\mu}_1} + \frac{1 - e^{-(\mu_2+\tilde{\mu}_2)T}}{\mu_2 + \tilde{\mu}_2}\Big).
\end{align*}
Similar expressions can be given for the covariances related to the
velocity and the cross terms. Here the eigenvalues of the nominal and
perturbed systems are (linear) functions of $\epsilon$ (and $\beta$)
that are related by the type of difference quotient chosen to
approximate the sensitivity.

In the case of a centered difference, $\mu_1 = \mu_1(\varepsilon)$ and
$\mu_2 = \mu_2(\varepsilon)$ are defined, in the obvious way, as
$\mu_1(\varepsilon) = (\beta+\varepsilon)/2 +
\sqrt{(\beta+\varepsilon)^2/4 - \omega^2}$
and
$\mu_2(\varepsilon) = (\beta+\varepsilon)/2 -
\sqrt{(\beta+\varepsilon)^2/4 - \omega^2}$
and hence $\tilde{\mu}_1 = \mu_1(-\varepsilon)$ and
$\tilde{\mu}_2 = \mu_2(-\varepsilon)$. In this case, we can write
$\pvar[X_T]= \psi(\mu_1(\varepsilon),\mu_2(\varepsilon))$ and
$\pvar[\tilde{X}_T] = \psi(\mu_1(-\varepsilon),\mu_2(-\varepsilon))$
where we define $\psi(\mu_1,\mu_2) = \phi(\mu_1,\mu_2,\mu_1,\mu_2)$.

The asymptotic dependence of
$\pvar[D(Z_t)] = \pvar[X_T] + \pvar[\tilde{X}_T] -2 \pcov[X_T,
\tilde{X}_T]$
on $\varepsilon$ can now be obtained by expanding the quantity of
interest in a series in $\varepsilon = 0$, using the representations
above. That is, for each terms appearing above we have 
\begin{equation*}\pvar[X_T] = \left. \psi \right|_{\varepsilon=0} +
  \left. \partial_\varepsilon \psi
  \right|_{\varepsilon=0} \varepsilon + \left. \partial^2_\varepsilon \psi
    \right|_{\varepsilon=0} \varepsilon^2 +
  O(\varepsilon^3),
\end{equation*}
\begin{equation*}\pvar[\tilde{X}_T] = \left. \psi
  \right|_{\varepsilon=0} - \left. \partial_\varepsilon \psi  \right|_{\varepsilon=0} \varepsilon +
 \left. \partial^2_\varepsilon \psi
    \right|_{\varepsilon=0} \varepsilon^2 + O(\varepsilon^3),
\end{equation*}
\begin{equation*}\pcov[X_T, \tilde{X}_T] = \left. \phi
  \right|_{\varepsilon=0} + \left. \partial_\varepsilon \phi \right|_{\varepsilon=0} \varepsilon +
  \left. \partial^2_\varepsilon \phi
    \right|_{\varepsilon=0} \varepsilon^2 + O(\varepsilon^3),
\end{equation*}
where $ \partial^k_\varepsilon$ denotes the $k$th derivative with
respect to $\varepsilon$. Noting that
$ \left. \psi \right|_{\varepsilon=0}$ is non-zero, it follows that
$\pvar[X_T - \tilde{X}_T] = O(1)$ for independently sampled paths. For
the common random path coupling, the zeroth order term in the
expansion for $\pvar[D(Z_t)]$ vanishes since
$ \left. \psi \right|_{\varepsilon=0} = \left. \phi
\right|_{\varepsilon=0}$.
In this particular case, the first order term,
$\left.\partial_\varepsilon \phi \right|_{\varepsilon=0} = 0$, also
vanishes since $\partial_{\mu_j}\phi = \partial_{\tilde{\mu_j}} \phi$
and $\tilde{\mu}_j^\prime(0) = -\mu_j^\prime(0)$, for $j = 1,2$.
Finally, noting that since $\partial_{\mu_j} \phi$ is not symmetric in
$\mu_j,\tilde{\mu}_j$, the second order term in the expansion for
$\pvar[D(Z_t)]$ does not vanish, yielding
$\pvar[D(Z_t)] = O(\varepsilon^2)$. Explicit expansions can also be
calculated for other standard first order differences and for the
other covariance terms with similar asymptotic rates observed, namely
$O(\varepsilon^2)$ for the common random path coupling and $O(1)$ for
independently sampled paths. 

\subsubsection{Time average observable}
\label{eg:time-average-observable}
\begin{figure}[h]
  \centering
  \includegraphics[width=\shrinkfac\columnwidth]{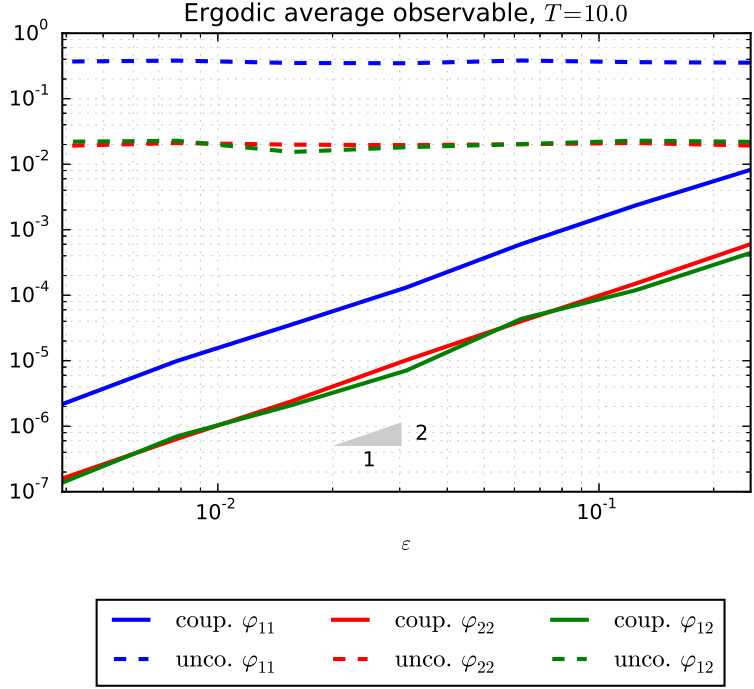}
  \caption{The sample variance of the components of $\pvar[D(Z_t)]$,
    for the time average observable $f(Z_t) = \overline{Z}$, are on
    the order of $O(\varepsilon^2)$ for optimally coupled paths, in
    contrast to $O(1)$ for independent paths. Here we consider of the
    sample variance of $\varphi_{11}$, $\varphi_{22}$, and
    $\varphi_{12}$, based on $M=10^3$ samples, for a central
    difference perturbation in $\beta$, at $\beta =1$, for time
    averages up to at $T = 10$.}
  \label{fig:langTimeAvgVariance}
\end{figure}

Let $\overline{X} = T^{-1} \int_0^T X_t \dd t$ and consider the
time average observable $f(Y_t) = (\overline{X}, \overline{V})$. As in
the case of the time average observable for the OU process, the
expectation can be exchanged with the integral in time, yielding
explicit expressions for the covariances as in \S
\ref{eg:finite-time-observable}. Investigations into the asymptotic
dependence of $\pvar[D(Z_t)]$ yield $O(\varepsilon^2)$ in the
optimally coupled case and $O(1)$ in the independent case. These rates
are observed experimentally in Figure~\ref{fig:langTimeAvgVariance}
where we consider $\pvar[D(Z_t)] = (\varphi)_{ij}$ (i.e.,
$\varphi_{11} = \pvar[\overline{X^1} - \overline{X^2}]$), based on
$M=10^3$ samples, for a central difference perturbation in $\beta$, at
$\beta =1$ (the underdamped case $\beta < 2\omega$). The time averages
are computed up to a final time $T = 10$ for sample paths from
Langevin dynamics, with fixed parameters $x_0=-1$, $v_0 = -0.1$,
$\omega=1$, $m=1$, and $\gamma = 1$, integrated using the BAOAB
method\cite{LeimkuhlerMatthews:2015md} with $\Delta t = 10^{-3}$.

\bibliography{uq_for_gle}
\end{document}